\documentclass{amsart}
\usepackage{amsfonts}
\usepackage{CJK,CJKnumb,CJKulem,times,dsfont,ifthen,mathrsfs,latexsym,amsfonts,color}
\usepackage{amsmath,amsthm,fontenc,amssymb,bm,graphicx,psfrag,listings,curves,extarrows}
\newtheorem{theorem}{Theorem}
\theoremstyle{plain}

\newtheorem{corollary}{Corollary}

\newtheorem{definition}{Definition}

\newtheorem{lemma}{Lemma}

\newtheorem{proposition}{Proposition}
\newtheorem{remark}{Remark}

\numberwithin{equation}{section}

\newcommand{\s}{\section}
\newcommand{\R}{\mathbb R}
\newcommand{\lab}{\label}
\newcommand{\bt}{\begin{theorem}}
\newcommand{\et}{\end{theorem}}
\newcommand{\bl}{\begin{lemma}}
\newcommand{\el}{\end{lemma}}
\newcommand{\bd}{\begin{definition}}
\newcommand{\ed}{\end{definition}}
\newcommand{\bc}{\begin{corollary}}
\newcommand{\ec}{\end{corollary}}
\newcommand{\bp}{\begin{proof}}
\newcommand{\ep}{\end{proof}}
\newcommand{\bo}{\begin{proposition}}
\newcommand{\eo}{\end{proposition}}
\newcommand{\br}{\begin{remark}}
\newcommand{\er}{\end{remark}}
\newcommand{\be}{\begin{equation}}
\newcommand{\ee}{\end{equation}}

\begin{document}
\title[Hardy-Sobolev critical exponents]{A nonlinear elliptic PDE with  multiple   Hardy-Sobolev critical exponents in $\R^N$}
\author{X. Zhong}
\address[X. Zhong]{Department of Mathematics, Sun Yat-sen University, Guangzhou 510275, China.}
\address[W. Zou]{Department of Mathematical Sciences, Tsinghua University, Beijing 100084, China.}
\email[X. Zhong]{zhongxuexiu1989@163.com}
\author{W. Zou}
\email[W. Zou]{wzou@math.tsinghua.edu.cn}
\thanks{Supported by NSFC (11025106, 11371212, 11271386) and the Both-Side Tsinghua Fund.}
\date{December 22, 2017}
\thanks{This paper is in final form and no version of it will be submitted
for publication elsewhere.}

\begin{abstract}
In this paper, we will study the following PDE in $\R^N$ involving  multiple   Hardy-Sobolev  critical exponents:
$$
\begin{cases}
\Delta u+\sum_{i=1}^{l}\lambda_i \frac{u^{2^*(s_i)-1}}{|x|^{s_i}}+u^{2^*-1}=0\;\hbox{in}\;\R^N,\\
u\in D_{0}^{1,2}(\R^N),
\end{cases}
$$
where $0<s_1<s_2<\cdots<s_l<2, 2^\ast:=\frac{2N}{N-2}, \; 2^\ast(s):=\frac{2(N-s)}{N-2}$ and  there exists some $k\in [1, l]$ such that $\lambda_i>0$ for $1\leq i\leq k$; $\lambda_i<0$ for $k+1\leq i\leq l$. We develop an interesting way to study this class of equations involving mixed sign parameters.
We  prove the existence and non-existence of the  positive ground state solution.   The regularity of the least-energy  solution are also  investigated.
\end{abstract}

\maketitle
\tableofcontents
\s{Introduction}
\renewcommand{\theequation}{1.\arabic{equation}}
\renewcommand{\theremark}{1.\arabic{remark}}
\renewcommand{\thedefinition}{1.\arabic{definition}}

\vskip0.1in
Consider the following problem:
\be\lab{PMS-0}
\begin{cases} \displaystyle
\Delta u+\sum_{i=1}^{l}\lambda_i \frac{u^{2^*(s_i)-1}}{|x|^{s_i}} +u^{2^*-1}=0\;\hbox{in}\;\R^N,\\
u(x)>0\;\;\hbox{in}\;\;\R^N,
\end{cases}
\ee
where $0<s_1<s_2<\cdots<s_l<2,  \lambda_i\in \R, 2^\ast:=\frac{2N}{N-2}, 2^\ast(s):=\frac{2(N-s)}{N-2}$.
We see that the nonlinearities involving multiple Hardy-Sobolev critical exponents and thus are not homogeneous.

Recall that on the half space $\R_+^N$,   Li and Lin consider the following problem  in \cite{LiLin.2012}:
\be\lab{poo}
\begin{cases}
\Delta u+\lambda \frac{u^{2^*(s)-1}}{|x|^{s_1}}+\frac{u^{2^*(s_2)-1}}{|x|^{s_2}}=0\quad\hbox{in}\;\R_+^N,\\
u(x)>0\;\;\;\hbox{in}\;\;\R_+^N, \;\;\;u(x)=0\;\;\hbox{on}\;\;\partial\R_+^N.
\end{cases}
\ee
 They show that (\ref{poo}) has a least-energy solution $u\in H_0^1(\R_+^N)$ provided that  $N\geq 3, 0<s_2<s_1<2, \lambda\in \R$.
 An earlier result for the special  case $s_2=0$ in  equation (\ref{poo}) is obtained   by Hsia, Lin and Wadade in \cite{HsiaLinWadade.2010}. Also   they study  the existence of  the least-energy solution.

\vskip0.1in

In the current paper, we consider the equation defined in the whole space $\R^N$ with multiple Hardy-Sobolev exponents.
It  seems that  the existence of least energy solution to (\ref{PMS-0}) is unknown.

\bt\lab{2014-12-12-mainth2}
Let $N\geq 3, 0<s_1<s_2<\cdots<s_l<2$. Suppose that  there exists some $k\in [1, l]$ such that $\lambda_i>0$ for $1\leq i\leq k$ and $\lambda_i<0$ for $k+1\leq i\leq l$. Furthermore, if $N=3$ and $k\neq l$, we assume that either $s_1<1$ or $1\leq s_1<2$ along with $\displaystyle\max\left\{|\lambda_{k+1}|,\cdots,|\lambda_l|\right\}$ small enough. Then the following problem
\be\lab{PMS}
\begin{cases}\displaystyle
\Delta u+\sum_{i=1}^{l}\lambda_i \frac{u^{2^*(s_i)-1}}{|x|^{s_i}} +u^{2^*-1}=0\;\hbox{in}\;\R^N,\\
u(x)>0\;\;\hbox{in}\;\;\R^N,
\end{cases}
\ee
has a least-energy solution  $u(x)$. Moreover, there exists a constant $C>0$ such that $u(x)\leq C(1+|x|^{2-N})$,
$|\nabla u(x)|\leq C(1+|x|^{1-N})$ and    $\displaystyle 0<\lim_{x\rightarrow 0}u(x)=\sup_{x\in \R^N} u(x)<\infty$.

\et

\vskip0.22in

\br  When all $\lambda_is$ are negative,  it is standard to prove that  there is  no least energy solution  to equation (\ref{PMS}). Thus, in the present paper, we always assume that $k\neq 0$. Besides, we may observe some different behaviors between $k=l$ and $k\neq l$. Actually, the result for $k=l$ can be proved in a much more direct way. However, it will encounter tough trouble. We prefer to adopt the way in present paper to study the case of $k=l$ since some results established in the process will be useful when we study the case of $k\neq l$.
\er

\br
Note that when $s_2>0$, we have $2^*(s_2)-1<2^*-1$. Then the subcritical equation
\be\lab{2015-3-17-ze4} \Delta \tilde{v}+\tilde{v}^{2^*(s_2)-1}=0\;\hbox{in}\;\Omega\ee
has no nontrivial solution if $\Omega=\R_+^N$. This result plays a  crucial role in \cite{LiLin.2012}.
 When we  consider the domain $\Omega=\R^N$, if $s_2>0$, we see that \eqref{2015-3-17-ze4} also has no nontrivial solution. For this case,   one can modify Li and Lin's arguments in \cite{LiLin.2012}  if $0<s_2<s_1<2$, and obtain the existence of ground state solution to problem $$ \Delta u+\lambda \frac{u^{2^*(s)-1}}{|x|^{s_1}}+\frac{u^{2^*(s_2)-1}}{|x|^{s_2}}=0, x\in \R^N, 0<u\in D_{0}^{1,2}(\R^N).$$   Note that this phenomenon will change  essentially when  $\Omega=\R^N$ and $s_2=0$. Since in this case, \eqref{2015-3-17-ze4} possesses a positive solution. Hence, when applying the blow-up method, ones  need a further detailed  arguments on the energy to deduce a contradiction. However, if we consider the problem \eqref{PMS} with $l>1$, i.e., the nonlinearities consist of multiple Hardy-Sobolev critical terms, the arguments of \cite{HsiaLinWadade.2010} can not be applied directly  to study the  equation (\ref{PMS}). Especially, when $k\neq l$, their arguments will fail.
\er


\br\lab{2015-3-26-r1}
The case of $k\neq l$ in  Theorem \ref{2014-12-12-mainth2} is much more complicated. It is not easy to exclude the``vanishing" phenomenon. We will apply a different method to study this case by considering  the following variant problem:
\be\lab{2015-3-28-we1}
    \begin{cases}
    \Delta u+u^{2^*-1}+\sum_{i=1}^{k}\lambda_i \frac{u^{2^*(s_i)-1}}{|x|^{s_i}}+\lambda\left(\sum_{i=k+1}^{l}|\lambda_i| \frac{u^{2^*(s_i)-1}}{|x|^{s_i}}\right)=0,\;x\in\;\R^N,\\
     u\in D_{0}^{1,2}(\R^N).
     \end{cases}
    \ee
    Denote
    \be
     D_k:=\left\{\mu\in \R\;:\;\hbox{problem \eqref{2015-3-28-we1} possesses a least energy solution when}\;\lambda=\mu\right\}.
    \ee
    We shall prove  $-1\in D_k$.
Basing on the results of section 2, we will apply the perturbation argument to deduce that $\emptyset\neq D_k$ is a set both open and closed.   Thus $D_k=\R$, and it follows that $-1\in D_k$,  which completes the proof.
\er


\br We remark that there are some works on the  Hardy-Sobolev critical elliptic equations with boundary singularities and on  the effect of curvature for  the best constant in the Hardy--Sobolev inequalities, see \cite{Ghou=1, Ghou=2, Ghou=3, HsiaLinWadade.2010, LiLin.2012,Lin-Wad}  and the references therein.   The limiting equations of \cite{Ghou=1, Ghou=2, Ghou=3, HsiaLinWadade.2010, LiLin.2012,Lin-Wad} are actually the form of
(\ref{PMS-0}) defined in a cone.

\er

This paper is organized as follows. In section 2, we will study the regularity of the nonnegative solution of \eqref{PMS}. In section 3 we will firstly study an  approximating problem of \eqref{PMS}. In section 4, we will introduce some interpolation inequalities and the Pohozaev identity for such equation. Finally, in section 5, we will prove the existence of ground state solution and  complete the proof of Theorem \ref{2014-12-12-mainth2}.


\vskip0.1in

\vskip0.1in

\s{The regularity of the  solution to  equation \eqref{PMS-0}}
\renewcommand{\theequation}{2.\arabic{equation}}
\renewcommand{\theremark}{2.\arabic{remark}}
\renewcommand{\thedefinition}{2.\arabic{definition}}
\bo\lab{2014-12-20-prop1}
Let $N\geq 3, 0<s_i<2, \lambda_i\neq 0, i=1,2,\cdots,l$ and set $s_0=0,\lambda_0=1$. Then any nonnegative $D_{0}^{1,2}(\R^N)-$solution $u$  of \be\lab{2014-12-20-e1}
\int_{\R^N}\nabla u\nabla \varphi dx=\sum_{i=0}^{l}\lambda_i\int_{\R^N}\frac{u^{2^*(s_i)-1}\varphi}{|x|^{s_i}}dx,\;\;\;\forall\;\varphi\in D_{0}^{1,2}(\R^N)
\ee
is of class $L_{loc}^{\infty}(\R^N)$.
\eo
\bp
Let $\chi$ be a cut-off function in a ball $\mathbb{B}_R(x_0)$. We take $\varphi=\chi^2uu_{M}^{2(t-1)}$, where $t>1,M>1$ and $u_M:=\min\{u,M\}$. Note that $\nabla u\nabla u_M=|\nabla u_M|^2$ and $\nabla u\nabla u_M u u_{M}^{2(t-1)-1}=|\nabla u_M|^2 u_{M}^{2(t-1)}$. Then by \eqref{2014-12-20-e1}, we have
\begin{align}\lab{2014-12-20-e2}
&\int_{\R^N}\nabla u\nabla \varphi dx\\
=&\int_{\R^N}\nabla u\cdot \left[2\chi\nabla \chi u u_{M}^{2(t-1)}+\chi^2 u_{M}^{2(t-1)}\nabla u+2(t-1)\chi^2 u u_{M}^{2(t-1)-1}\nabla u_M\right]\nonumber\\
=&\int_{\R^N}\chi^2|\nabla u|^2 u_{M}^{2(t-1)}dx+2(t-1)\int_{\R^N}\chi^2|\nabla u_M|^2u_{M}^{2(t-1)}dx\nonumber\\
&+2\int_{\R^N}\chi u u_{M}^{2(t-1)}\nabla \chi \nabla u dx\nonumber\\
=&\sum_{i=0}^{l}\lambda_i\int_{\R^N}\frac{u^{2^*(s_i)}}{|x|^{s_i}}\chi^2 u_{M}^{2(t-1)}dx.\nonumber
\end{align}
By the Young's inequality, we have
\be\lab{2014-12-20-e3}
\left|(\chi\nabla u) \cdot (u\nabla \chi)\right|\leq \left|\nabla \chi\right|^2u^2+\frac{1}{4}\chi^2\left|\nabla u\right|^2.
\ee
Thus,
\begin{align}\lab{2014-12-20-e4}
&\int_{\R^N}\chi^2|\nabla u|^2 u_{M}^{2(t-1)}dx+2(t-1)\int_{\R^N}\chi^2|\nabla u_M|^2u_{M}^{2(t-1)}dx\\
& \leq \sum_{i=0}^{l}\lambda_i\int_{\R^N}\frac{u^{2^*(s_i)}}{|x|^{s_i}}\chi^2 u_{M}^{2(t-1)}dx+2\int_{\R^N}u_{M}^{2(t-1)}\left[|\nabla \chi|^2u^2+\frac{1}{4}\chi^2|\nabla u|^2\right]dx\nonumber
\end{align}
and it follows that
\begin{align}\lab{2014-12-20-e5}
&\int_{\R^N}\chi^2|\nabla u|^2 u_{M}^{2(t-1)}dx+4(t-1)\int_{\R^N}\chi^2|\nabla u_M|^2u_{M}^{2(t-1)}dx\\
&\leq  2\sum_{i=0}^{l}\lambda_i\int_{\R^N}\frac{u^{2^*(s_i)}}{|x|^{s_i}}\chi^2 u_{M}^{2(t-1)}dx+4\int_{\R^N}u_{M}^{2(t-1)}|\nabla \chi|^2u^2 dx.\nonumber
\end{align}
Now, we take $t=\frac{2^*(\overline{s})}{2}>1$ with $\overline{s}:=\max\{s_0,s_1,\cdots, s_l\}=s_l<2$ for simplicity. Consider $w_M:=\chi u u_{M}^{t-1}$, by the H\"older inequality and the Hardy-Sobolev inequality, formula \eqref{2014-12-20-e5} yields that
\begin{align}\lab{2014-12-20-e6}
&\left(\int_{\R^N}\frac{1}{|x|^{s_i}}(w_M)^{2^*(\overline{s})}dx\right)^{\frac{2}{2^*(\overline{s})}}
\leq C\left(\int_{\R^N}\frac{1}{|x|^{s_i}}(w_M)^{2^*(s_i)}dx\right)^{\frac{2}{2^*(s_i)}}
 \leq C\int_{\R^N}|\nabla w_M|^2 dx\\
& \leq C_1 \Big[\int_{\R^N}|\nabla \chi|^2u^2u_{M}^{2(t-1)}dx+\int_{\R^N}\chi^2|\nabla u|^2u_{M}^{2(t-1)}dx+
 (t-1)^2\int_{\R^N}\chi^2u_{M}^{2(t-1)}|\nabla u_M|^2dx\Big]\nonumber\\
& \leq C_2 t\Big[\sum_{i=0}^{l}|\lambda_i|\int_{\R^N}\frac{u^{2^*(s_i)}}{|x|^{s_i}}\chi^2u_{M}^{2(t-1)}dx+  \int_{\R^N}u_{M}^{2(t-1)}|\nabla\chi|^2u^2dx\Big].\nonumber
\end{align}
By the H\"older inequality,
\begin{align}\lab{2014-12-20-e7}
&\left|\lambda_i \int_{\R^N}\frac{u^{2^*(s_i)}}{|x|^{s_i}}\chi^2u_{M}^{2(t-1)}dx \right|\\
\leq& |\lambda_i|\left[\int_{\mathbb{B}_{R_0}(x_0)}\frac{u^{2^*(s_i)}}{|x|^{s_i}}dx\right]^{\frac{2^*(s_i)-2}{2^*(s_i)}}
\cdot \left[\int_{\R^N}\frac{1}{|x|^{s_i}}\left(\chi u u_{M}^{t-1}\right)^{2^*(s_i)}dx\right]^{\frac{2}{2^*(s_i)}},\nonumber
\end{align}
then by the absolute continuity of the integral, we see that there exists some $R_0>0$ small enough such that
\be\lab{2014-12-20-e8}
\frac{2^*(\overline{s})}{2}C_2 |\lambda_i| \left[\int_{\mathbb{B}_{R_0}(x_0)}\frac{u^{2^*(s_i)}}{|x|^{s_i}}dx\right]^{\frac{2^*(s_i)-2}{2^*(s_i)}}<\frac{1}{l+2}\;\hbox{for all}\;i=0,1,2,\cdots,l.
\ee
Hence, for such $R_0$, we have
\begin{align}\lab{2014-12-20-e9}
&\left(\int_{\R^N}\frac{1}{|x|^{s_i}}(w_M)^{2^*(s_i)}dx\right)^{\frac{2}{2^*(s_i)}}\\
\leq&\frac{1}{l+2}\sum_{j=0}^{l}\left(\int_{\R^N}\frac{1}{|x|^{s_j}}(w_M)^{2^*(s_j)}dx\right)^{\frac{2}{2^*(s_j)}}+\tilde{C}_2\int_{\R^N}|\nabla \chi|^2 u^{2^*(\overline{s})}dx\;\nonumber
\end{align}
for all  $\;i=0,1,\cdots,l.$ It follows that
\be\lab{2014-12-20-e10}
\sum_{i=0}^{l}\left(\int_{\R^N}\frac{1}{|x|^{s_i}}(w_M)^{2^*(s_i)}dx\right)^{\frac{2}{2^*(s_i)}}
\leq C_3\int_{\R^N}|\nabla \chi|^2 u^{2^*(\overline{s})}dx.
\ee
Let $M$ go to infinity, we obtain that
\be\lab{2014-12-20-e11}
u\in L^{\frac{2^*(\overline{s})}{2} 2^*(s_i)}\left(\mathbb{B}_{\frac{R}{2}}(x_0), \frac{dx}{|x|^{s_i}}\right), \;\;i=0,1,\cdots,l.
\ee
By the arbitrariness  of $x_0$, we obtain that
\be\lab{2014-12-20-e12}
u\in L_{loc}^{\frac{2^*(\overline{s})}{2} 2^*(s_i)}\left(\R^N, \frac{dx}{|x|^{s_i}}\right), \;\;i=0,1,\cdots,l.
\ee
Now, for any $R>0, 0<r<1$, we take a cut-off function $0<\chi\leq 1$ in $\mathbb{B}_{R+r}$ such that $\chi\equiv 1$ in $\mathbb{B}_R$ and $|\nabla \chi|\leq \frac{2}{r}$ in $\mathbb{B}_{R+r}$. Set
\be\lab{2014-12-20-e14}
\sigma_i:=\frac{2^*(\overline{s}) 2^*(s_i)}{2[2^*(s_i)-2]},\quad i=0,1,\cdots,l.
\ee
We note that
\be\lab{2014-12-20-bue1}
\frac{2^*(s_i)(\sigma_i-1)}{2\sigma_i}>1\;\hbox{for all}\;i=0,1,\cdots,l,
\ee
we can take proper constants $q_i\leq 2^*(s_i)$ such that
\be\lab{2014-12-20-bue2}
\frac{q_i(\sigma_i-1)}{2\sigma_i}\equiv const>1.
\ee
By the  H\"older inequality, we have
\begin{align}\lab{2014-12-20-e15}
&\int_{\R^N}\frac{1}{|x|^{s_i}} u^{2^*(s_i)-2}\chi^2 u_{M}^{2t} dx\\
\leq&\left[\int_{\mathbb{B}_{R+r}}\frac{1}{|x|^{s_i}}u^{[2^*(s_i)-2]\sigma_i}\chi^2 dx\right]^{\frac{1}{\sigma_i}} \left[\int_{\mathbb{B}_{R+r}}\frac{1}{|x|^{s_i}}\chi^2 u_{M}^{\frac{2t\sigma_i}{\sigma_i-1}}dx\right]^{\frac{\sigma_i-1}{\sigma_i}}\nonumber\\
\leq&C_4\left[\int_{\mathbb{B}_{R+r}}\frac{1}{|x|^{s_i}}\chi^2 u^{\frac{2t\sigma_i}{\sigma_i-1}}dx\right]^{\frac{\sigma_i-1}{\sigma_i}},\nonumber
\end{align}
provided that $\displaystyle u\in L^{\frac{2t\sigma_i}{\sigma_i-1}}\left(\mathbb{B}_{R+r}, \frac{dx}{|x|^{s_i}}\right)$.
Here we remark that by the H\"oler inequality, $C_4$ should depend on the volume of the ball $\mathbb{B}_{R+r}$. However, since $r<1$, we can choose some suitable $C_4$ that independent of $r$.
Noting that the right hand side of \eqref{2014-12-20-e15} is independent of $M$, by letting $M$ go to infinity, we indeed obtain that
\begin{align}\lab{2014-12-20-e16}
&\int_{\R^N}\frac{1}{|x|^{s_i}}u^{2^*(s_i)}\chi^2 u_{M}^{2(t-1)}dx
\leq\int_{\R^N}\frac{1}{|x|^{s_i}}u^{2^*(s_i)}\chi^2 u^{2(t-1)}dx\\
\leq&C_4\left[\int_{\mathbb{B}_{R+r}}\frac{1}{|x|^{s_i}} u^{\frac{2t\sigma_i}{\sigma_i-1}}dx\right]^{\frac{\sigma_i-1}{\sigma_i}}.\nonumber
\end{align}
On the other hand,
\begin{align}\lab{2014-12-20-e17}
&\int_{\R^N}|\nabla \chi|^2u^2u_{M}^{2(t-1)}dx\\
\leq&\left(\frac{2}{r}\right)^2\left(\int_{\mathbb{B}_{R+r}}\frac{1}{|x|^{s_i}}u^{\frac{2t\sigma_i}{\sigma_i-1}}dx\right)^{\frac{\sigma-1}{\sigma}} \left(\int_{\mathbb{B}_{R+r}}|x|^{\eta_i}dx\right)^{\frac{1}{\sigma_i}},\nonumber
\end{align}
where
$$\eta_i=s_i(\sigma_i-1)\geq 0.$$
Hence,
\be\lab{2014-12-20-e18}
\int_{\R^N}|\nabla \chi|^2u^2u_{M}^{2(t-1)}dx\leq C_5 r^{-2}\left(\int_{\mathbb{B}_{R+r}}\frac{1}{|x|^{s_i}} u^{\frac{2t\sigma_i}{\sigma_i-1}}dx\right)^{\frac{\sigma_i-1}{\sigma_i}}.
\ee
Recalling that $q_i\leq 2^*(s_i)$ and $\chi\equiv 1$ in $\mathbb{B}_R$, by the H\"older inequality,
inserting \eqref{2014-12-20-e16},\eqref{2014-12-20-e18} into \eqref{2014-12-20-e6} and then letting $M$ go to infinity, we obtain that
\be\lab{2014-12-20-e20}
\sum_{i=0}^{l}\left(\int_{\mathbb{B}_R} \frac{1}{|x|^{s_i}}u^{q_it}dx\right)^{\frac{1}{q_it}}
\leq C_{6}^{\frac{1}{t}} t^{\frac{1}{2t}} r^{-\frac{1}{t}}\sum_{i=0}^{l}\left(\int_{\mathbb{B}_{R+r}}\frac{1}{|x|^{s_i}} u^{tq_{i,0}}dx\right)^{\frac{1}{tq_{i,0}}},
\ee
where by \eqref{2014-12-20-bue2},
$$q_{i,0}:=\frac{2\sigma_i}{\sigma_i-1}<q_i.$$
Recalling \eqref{2014-12-20-bue2} again, we have that $\tau:=\frac{q_i}{q_{i,0}}>1$ is independent of $i$. Define $t=\tau^j, R=1$ and $r_j=2^{-j}, j\geq 1$, applying iteration, \eqref{2014-12-20-e20} yields
\be\lab{2014-12-20-e21}
\sum_{i=0}^{l}\left[\int_{\mathbb{B}_{1+2^{-j+1}}} \frac{1}{|x|^{s_i}}u^{q_i\tau^{j+1}}dx\right]^{\frac{1}{q_i\tau^{j+1}}}
\leq \prod_{k=0}^{j+1}(C_6\tau^{\frac{k}{2}}2^k)^{\tau^{-k}}\sum_{i=0}^{l}\left[\int_{\mathbb{B}_{\frac{3}{2}}}\frac{1}{|x|^{s_i}} u^{q_{i}}dx\right]^{\frac{1}{q_{i}}}.
\ee
Denote
$$\Theta:=\prod_{k=0}^{\infty}(C_6\tau^{\frac{k}{2}}2^k)^{\tau^{-k}},$$
we have
\begin{align}\lab{2014-12-20-e22}
\ln \Theta=\ln C_6\sum_{k=0}^{\infty}\frac{1}{\tau^k}+\big(\ln 2+\frac{1}{2}\ln \tau\big)\sum_{k=0}^{\infty}\frac{k}{\tau^k}.
\end{align}
It is easy to see $\Theta<\infty$ due to the fact of $\tau>1$.
Hence, letting $j$ go to infinity in \eqref{2014-12-20-e21}, noting that $s_i\geq 0$, we obtain that
\be\lab{2014-12-20-e23}
\sup_{\mathbb{B}_1}\;u\leq \frac{1}{l+1}\Theta \sum_{i=0}^{l}\left[\int_{\mathbb{B}_{\frac{3}{2}}}\frac{1}{|x|^{s_i}} u^{q_{i}}dx\right]^{\frac{1}{q_{i}}}.
\ee
\ep
Then we have the following result.

\bl\lab{2014-12-18-wl1}
Let $N\geq 3, 0<s_i<2, \lambda_i>0,i=1,2,\cdots, l$. Then any nonnegative solution of
\be\lab{2014-12-18-we1}
\Delta u+\sum_{i=1}^{l}\lambda_i \frac{u^{2^*(s_i)-1}}{|x|^{s_i}} +u^{2^*-1}=0\;\hbox{in}\;\R^N\backslash\{0\}
\ee
satisfying
\be\lab{2014-12-18-we2}
0<\liminf_{x\rightarrow 0} u(x)\leq \limsup_{x\rightarrow 0}u(x)<+\infty.
\ee
\el
\bp
By the standard elliptic estimation, we have that $u\in C^\infty(\R^N)\backslash \{0\}$. Then by \cite[Lemma 4.2]{ChouChu.1993}, take some $r>0$,  we see that $t\mapsto \min_{|x|=t}u(x)$ is concave in $t^{2-N}$ for $t\in(0,r)$. Hence,
\be\lab{2014-12-18-we3}
u(x)\geq \min_{|x|=r}u(x)>0\;\hbox{for all}\;x\in \overline{B_r}\backslash\{0\},
\ee
and thus
\be\lab{2014-12-18-we4}
\liminf_{x\rightarrow 0}u(x)\geq \min_{|x|=r}u(x)>0.
\ee
On the other hand, by Proposition \ref{2014-12-20-prop1}, $u(x)$ is of class $L_{loc}^{\infty}(\R^N)$. Hence, the proof of this lemma is completed.
\ep

\bl\lab{2014-12-20-wl1}
Let $N\geq 3, 0<s_i<2, \lambda_i>0, i=1,2,\cdots,l$ and set $s_0=0,\lambda_0=1$. Then any nonnegative $D_{0}^{1,2}(\R^N)-$solution  of \eqref{2014-12-20-e1} satisfying
$$0<\liminf_{|x|\rightarrow \infty} |x|^{N-2}u(x)\leq \limsup_{|x|\rightarrow \infty}|x|^{N-2}u(x)<\infty,$$
i.e.,
$u=O(\frac{1}{|x|^{N-2}})$ when $|x|\rightarrow +\infty$.
\el
\bp
When $u$ is a nonnegative solution of \eqref{2014-12-20-e1}, a direct computation shows that its Kelvin Transform $\displaystyle v(x):=\left|x\right|^{-(N-2)}u\left(\frac{x}{|x|^2}\right)$ is also a nonnegative solution of
\eqref{2014-12-20-e1}. Then by Lemma \ref{2014-12-18-wl1},  we have
$$0<\liminf_{x\rightarrow 0} v(x)\leq \limsup_{x\rightarrow 0}v(x)<+\infty,$$
which implies the results of this Lemma.
\ep

\br\lab{2014-12-27-zr1}
Indeed, even for the case of $\lambda_i<0$, if $u(x)$ is a nonnegative solution of \eqref{2014-12-20-e1}, the corresponding Kelvin Transform $\displaystyle v(x):=|x|^{-(N-2)}u\left(\frac{x}{|x|^2}\right)$ is also a nonnegative solution. Then by Proposition \ref{2014-12-20-prop1}, $|v(x)|\leq C$ for $|x|<1$. Thus,
\be\lab{2014-12-17-zbue1}
|u(x)|\leq C |x|^{2-N}\quad\;\hbox{for}\;|x|\geq 1.
\ee
Furtherore, we can obtain that
\be\lab{2014-12-17-zbue2}
|u(x)|\leq C \left(1+|x|^{2-N}\right)\;\hbox{for}\;x\in \R^N.
\ee
And the standard gradient estimate implies that
\be\lab{2014-12-17-zbue3}
|\nabla u(x)|\leq C|x|^{1-N}\;\hbox{for}\;|x|\geq 1.
\ee
\er
\s{Approximating problems}
\renewcommand{\theequation}{3.\arabic{equation}}
\renewcommand{\theremark}{3.\arabic{remark}}
\renewcommand{\thedefinition}{3.\arabic{definition}}
Assume that $0=s_0<s_1<s_2<\cdots<s_l<2$.
 Let $0<\varepsilon<s_1$ and define
\be\lab{2014-12-12-we1}
a_{i,\varepsilon}(x):=\begin{cases}
\frac{1}{|x|^{s_i-\varepsilon}}\quad &\hbox{for}\;|x|<1,\\
\frac{1}{|x|^{s_i+\varepsilon}}&\hbox{for}\;|x|\geq 1,
\end{cases}\quad i=1,2,\cdots,l.
\ee
We also denote $a_{i,0}(x)=\frac{1}{|x|^{s_i}}$. Then it is easy to see that $a_{i,\varepsilon}(x)$ is decreasing with respect  to $\varepsilon\in [0,s_1)$.

\bl\lab{2014-12-12-wl1}
Let $0\leq \varepsilon<s_1$, then
for any $\displaystyle u\in D_{0}^{1,2}(\R^N)$ and $\displaystyle i\in \{1,2,\cdots,l\},$
 $\displaystyle\int_{\R^N} a_{i,\varepsilon}(x)|u|^{2^*(s_i)}dx$ is well defined and decreasing by $\varepsilon$.
\el
\bp
See  \cite[Lemma 14]{ZhongZou.2017}.
\ep

Denote by $L^{p}\left(\R^N,a_{i,\varepsilon}(x)dx\right)$ the space of $L^{p}$-integrable functions with respect to the measure $a_{i,\varepsilon}(x)dx$ and the corresponding norm is indicated by $$\left|u\right|_{p,i, \varepsilon}:=\ \left(\int_{\R^N} a_{i,\varepsilon}(x)|u|^pdx\right)^{\frac{1}{p}},\quad p>1.$$
Then we have the following result on the compact embedding.
\bl\lab{2014-12-3-l2}
For any $\varepsilon\in (0, s_1)$, the embedding $D_{0}^{1,2}(\R^N)\hookrightarrow L^{2^*(s_i)}\left(\R^N,a_{i,\varepsilon}(x)dx\right)$ is compact.
\el
\bp We refer to \cite[Lemma 16]{ZhongZou.2017} for the details.
\ep

We note that for any compact set $\Omega\subset \R^N$ with  $0\not\in\bar{\Omega}$, we have that $a_{i,\varepsilon}(x)\rightarrow a_{i,0}(x)$ uniformly for $i\in\{1,2,\cdots, l\}$ and $x\in\Omega$ as $\varepsilon\rightarrow 0$. Now, for any $0<\varepsilon< s_1$ fixed, let us consider the following problem:
\be\lab{Pva}
\begin{cases}\displaystyle
\Delta u+\sum_{i=1}^{l}\lambda_i  a_{i,\varepsilon}(x)u^{2^*(s_i)-1} +u^{2^*-1}=0\;\hbox{in}\;\R^N\backslash\{0\},\\
u(x)>0\;\;\hbox{in}\;\;\R^N\backslash \{0\},\quad u\in D_{0}^{1,2}(\R^N),
\end{cases}
\ee
whose energy functional is given by
\be
\Phi_\varepsilon(u)=\frac{1}{2}\int_{\R^N}|\nabla u|^2 dx-\sum_{i=1}^{\lambda_i}\frac{1}{2^*(s_i)} \int_{\R^N}a_{i,\varepsilon}(x)|u|^{2^*(s_i)}dx-\frac{1}{2^*} \int_{\R^N}|u|^{2^*}dx.
\ee

\subsection{Nehari Manifold $\mathcal{N}_\varepsilon$}
Consider the corresponding Nehari manifold
\be\lab{2014-12-16-e1}
\mathcal{N}_\varepsilon:=\left\{u\in D_{0}^{1,2}(\R^N)\backslash \{0\}\;: J_\varepsilon(u)=0\right\},
\hbox{ where }
J_\varepsilon(u):=\langle \Phi'_\varepsilon(u), u \rangle.
\ee
The following properties of $\mathcal{N}_\varepsilon$ are basic and the proofs are standard. For the reader's convenience, we  give the details.
\bl\lab{2014-12-16-l1}
Let $N\geq 3, 0=s_0<s_1<s_2<\cdots<s_l<2, \lambda_0=1$, there exists some $1\leq k\leq l$ such that $\lambda_i>0$ for $1\leq i\leq k$ and $\lambda_i<0$ for $k+1\leq i\leq l$. Then for any $u\in D_{0}^{1,2}(\R^N)\backslash \{0\}$, there exists a unique $t_{(\varepsilon,u)}>0$ such that $t_{(\varepsilon,u)}u\in \mathcal{N}_\varepsilon$. Further, $\mathcal{N}_\varepsilon$ is closed and bounded away from $0$. Moreover, if $k= l$, for any fixed $u\in D_{0}^{1,2}(\R^N)\backslash \{0\}$, $t=t_{(\varepsilon,u)}$ is strictly increasing by $\varepsilon$ in $[0, s_1)$.
\el
\bp Firstly, we consider the case of that $k=l$.
For any $0\neq u\in D_{0}^{1,2}(\R^N)$, we set
\be\lab{2014-12-16-e3}
b_{i,\varepsilon}(u):=\lambda_i\int_{\R^N}a_{i,\varepsilon}(x)|u|^{2^*(s_i)}dx>0,\quad i=1,2,\cdots,l.
\ee
We note that $b_{i,\varepsilon}(u)$ is strictly decreasing by $\varepsilon$ due to the monotonicity of $a_{i,\varepsilon}(x)$.
Since
\be\lab{2014-12-16-e4}
\Phi_\varepsilon(tu)=\frac{1}{2}\|u\|^2 t^2-\sum_{i=1}^{l}\frac{b_{i,\varepsilon}(u)}{2^*(s_i)}t^{2^*(s_i)}-\frac{1}{2^*}|u|_{2^*}^{2^*} t^{2^*},
\ee
by a direct computation, we  see that $\frac{d\Phi_{\varepsilon}(tu)}{dt}=0$ has a unique solution $t_{(\varepsilon,u)}>0$. Precisely, $t_{(\varepsilon,u)}$ is implicitly  given by the following algebraic equation
\be\lab{2014-12-16-e5}
\|u\|^2-\sum_{i=1}^{l} b_{i,\varepsilon}(u)t^{2^*(s_i)-2}-|u|_{2^*}^{2^*}t^{2^*-2}=0.
\ee
By Sobolev inequality, it is easy to see that there exists $\delta_\varepsilon>0$ such that $t_{(\varepsilon,u)}\geq \delta_\varepsilon$ for any $u$ satisfying $\|u\|=1$. Hence, $\mathcal{N}_\varepsilon$ is bounded away form $0$. Now, we prove that $t=t_{(\varepsilon,u)}$ is increasing by $\varepsilon$.
Assume that $0\leq \varepsilon_1<\varepsilon_2<s_1$, then we see that there exists a unique $t_1$ and $t_2$ such that
\be\lab{2014-12-16-e6}
J_{\varepsilon_1} (t_1u)=J_{\varepsilon_2} (t_2u)=0.
\ee
Recalling that  $b_{i,\varepsilon}(u)$ is strictly decreasing by $\varepsilon$, we have that
\be\lab{2014-12-16-e7}
J_{\varepsilon_2}(t_1u)>J_{\varepsilon_1}(t_1u)=0=J_{\varepsilon_2}(t_2u).
\ee
Noting that $\displaystyle J_{\varepsilon_2}(tu)\rightarrow -\infty$ as $t\rightarrow +\infty$,
then by the uniqueness of $t_{\varepsilon_2, u}$, we see that $t_{(\varepsilon_2, u)}=t_2>t_1$. Hence, $t_{(\varepsilon, u)}$ is strictly increasing by $\varepsilon$ in $[0,s_1)$.

\vskip 0.2in
\noindent
   Secondly, we consider the case of that $k\neq l$. For the convenience, we denote $b_{0,\varepsilon}(u)\equiv |u|_{2^*}^{2^*}$. Then we see that
\be\lab{2014-12-25-e1}
\Phi_{\varepsilon}(tu)=\frac{1}{2}\|u\|^2 t^2-\sum_{i=0}^{l}\frac{1}{2^*(s_i)}b_{i,\varepsilon}(u)t^{2^*(s_i)}.
\ee
For a given $u\neq 0$, we have
\be\lab{2014-12-25-e2}
\frac{d}{dt}\Phi_\varepsilon(tu):=tf(t),
\ee
where
\be\lab{2014-12-25-e3}
f(t):=\|u\|^2-\sum_{i=0}^{l}b_{i,\varepsilon}(u)t^{2^*(s_i)-2}.
\ee
Noting that $\frac{d}{dt}\Phi_\varepsilon(tu)=0$ with $t>0$ if and only if $f(t)=0$, hence the existence of $t_{(\varepsilon,u)}$ follows  easily from the continuity of $f(t)$  and the facts  that $\displaystyle f(0)=\|u\|^2>0, \lim_{t\rightarrow +\infty}f(t)=-\infty$.

Now, we shall prove the uniqueness of $t_{(\varepsilon,u)}$. Set $\mathcal{A}:=\{t>0\;:\; f(t)=0\}$. Then, we see that $\mathcal{A}\neq \emptyset$. Let $t_0:=\inf\; \mathcal{A}$, then by it is easy to see that $t_0>0$ and $t_0\in \mathcal{A}$, i.e., $t_0$ is the minimal positive root of $f(t)=0$. Hence, we have $f(t)>0$ for $0<t<t_0$ and $f(t_0)=0$. Next, we will show that $f'(t)<0$ for $t>t_0$ and thus $f(t)<f(t_0)=0$ for $t>t_0$. Indeed, $$\displaystyle f'(t)=-\sum_{i=0}^{l}[2^*(s_i)-2]b_{i,\varepsilon}
(u)t^{2^*(s_i)-3}:=-t^{2^*(s_{k+1})-3}
\left[\sum_{i=0}^{l}[2^*(s_i)-2]
b_{i,\varepsilon}(u)t^{2^*(s_i)-2^*(s_{k+1})}\right]$$
and thus we only need to prove that
\be\lab{2014-12-25-e4}
g(t):=\sum_{i=0}^{l}[2^*(s_i)-2]
b_{i,\varepsilon}(u)t^{2^*(s_i)-2^*(s_{k+1})}>0\;\hbox{for}\;t>t_0.
\ee
By $f(t_0)=0$, we have
\be\lab{2014-12-25-e5}
b_{k+1,\varepsilon}(u)=\|u\|^2 t_{0}^{2-2^*(s_{k+1})}-\sum_{i\neq k+1}b_{i,\varepsilon}(u)t_{0}^{2^*(s_i)-2^*(s_{k+1})}.
\ee
Noting that $b_{i,\varepsilon}(u)>0, 2^*(s_i)>2^*(s_{k+1})$ for $0\leq i\leq k$, then $b_{i,\varepsilon}(u)t^{2^*(s_i)-2^*(s_{k+1})}$ is increasing by $t$ in $(t_0,+\infty)$. We also note that $b_{i,\varepsilon}(u)<0, 2^*(s_i)<2^*(s_{k+1})$ for  $k+1<i\leq l$. Hence, for this case, we also obtain that $\displaystyle b_{i,\varepsilon}(u)t^{2^*(s_i)-2^*(s_{k+1})}$ is increasing by $t$ in $(t_0,+\infty)$. It is trivial that  $2^*(s_i)-2>0$ for $i=0,1,\cdots, l$. Hence,
\begin{align}\lab{2014-12-25-e6}
g(t):=&\sum_{i\neq k+1}[2^*(s_i)-2]
b_{i,\varepsilon}(u)t^{2^*(s_i)-2^*(s_{k+1})}+[2^*(s_{k+1})-2]b_{k+1,\varepsilon}(u)\\
>&\sum_{i\neq k+1}[2^*(s_i)-2]
b_{i,\varepsilon}(u)t_{0}^{2^*(s_i)-2^*(s_{k+1})}+[2^*(s_{k+1})-2]b_{k+1,\varepsilon}(u)\nonumber\\
=&\sum_{i\neq k+1}[2^*(s_i)-2]
b_{i,\varepsilon}(u)t_{0}^{2^*(s_i)-2^*(s_{k+1})}\nonumber\\
&+[2^*(s_{k+1})-2]\left[\|u\|^2 t_{0}^{2-2^*(s_{k+1})}-\sum_{i\neq k+1}b_{i,\varepsilon}(u)t_{0}^{2^*(s_i)-2^*(s_{k+1})}\right]\nonumber\\
=&[2^*(s_{k+1})-2]\|u\|^2 t_{0}^{2-2^*(s_{k+1})}+\nonumber\\
&\quad \sum_{i\neq k+1}[2^*(s_i)-2^*(s_{k+1})]
b_{i,\varepsilon}(u)t_{0}^{2^*(s_i)-2^*(s_{k+1})}.\nonumber
\end{align}
Noting  $2^*(s_{k+1})-2>0$ and   $\displaystyle [2^*(s_i)-2^*(s_{k+1})]
b_{i,\varepsilon}(u)>0$ for $i\neq k+1$.
Thus, \eqref{2014-12-25-e4} is proved and thereby we obtain the uniqueness of $t_{(\varepsilon,u)}$.
However,   for the case of  $k\neq l$,  we can not obtain the monotonicity of $t_{(\varepsilon,u)}$ by $\varepsilon$ in $[0,s_1)$.

\ep

\bl\lab{2014-12-16-l2}
Under the assumption of Lemma \ref{2014-12-16-l1},    any $(PS)_c$ sequence of $\Phi_\varepsilon(u)$
is bounded in $D_{0}^{1,2}(\R^N)$.
\el
\bp
Since $\{u_n\}$ is a $(PS)_c$ sequence, i.e.,
$\displaystyle
\Phi_\varepsilon(u_n)=c+o(1)$ and
$\langle \Phi'_\varepsilon(u_n),u_n \rangle=o(1)\|u_n\|,$ we have
\begin{align}\lab{2014-12-16-e11}
&c+o(1)+o(1)\|u_n\|\\
=&\left[\frac{1}{2}-\frac{1}{2^*(s_k)}\right]\|u_n\|^2+\sum_{i=1}^{l}\left[\frac{1}{2^*(s_k)}-\frac{1}{2^*(s_i)}\right]
b_{i,\varepsilon}(u_n)+\left[\frac{1}{2^*(s_k)}-\frac{1}{2^*}\right]|u_n|_{2^*}^{2^*}\nonumber\\
>&\left[\frac{1}{2}-\frac{1}{2^*(s_k)}\right]\|u_n\|^2,\nonumber
\end{align}
which implies that $\{u_n\}$ is bounded in $D_{0}^{1,2}(\R^N)$.
\ep

Define
\be\lab{2014-12-5-e3}
c_\varepsilon:=\inf_{u\in \mathcal{N}_\varepsilon}\Phi_\varepsilon(u)
\ee
and
\be\lab{2014-12-5-e4}
\delta_\varepsilon:=\inf_{u\in \mathcal{N}_\varepsilon}\|u\|.
\ee
\bl\lab{2014-12-16-l3}
Under the assumptions of Lemma \ref{2014-12-16-l1} and furthermore we suppose  that  $k=l$, i.e., all $\lambda_is$ are positive. Then $\delta_\varepsilon$ is strictly increasing by $\varepsilon\in [0,s_1)$, i.e., $0<\delta_0\leq \delta_{\varepsilon_1}< \delta_{\varepsilon_2}$ provided $0\leq \varepsilon_1<\varepsilon_2<s_1$.
\el
\bp
It follows from  the strictly monotonicity of $t_{(\varepsilon, u)}$ in   Lemma \ref{2014-12-16-l1}.
\ep

\bl\lab{2014-12-16-l4}
Under the assumptions of Lemma \ref{2014-12-16-l1}. Let $\{u_n\}$ be a $(PS)_c$ sequence of $\Phi_\varepsilon\big|_{\mathcal{N}_\varepsilon}$ i.e.,
$$\begin{cases}\Phi_\varepsilon(u_n)\rightarrow c\\ \Phi'_\varepsilon\big|_{\mathcal{N}_\varepsilon}(u_n)\rightarrow 0\;\hbox{in the dual space of}\;D_{0}^{1,2}(\R^N) \end{cases},$$
then $\{u_n\}$ is also a $(PS)_c$ sequence of $\Phi_\varepsilon$.
\el
\bp
For any $u\in \mathcal{N}_\varepsilon$, we have
\be\lab{2014-12-16-e12}
J_\varepsilon(u)=\|u\|^2-\sum_{i=1}^{l} b_{i,\varepsilon}(u)-|u|_{2^*}^{2^*}=0.
\ee
Consider the case of $k=l$, we have
\begin{align}\lab{2014-12-16-e13}
\langle J'_\varepsilon(u), u\rangle =&2\|u\|^2-\sum_{i=1}^{l} 2^*(s_i)b_{i,\varepsilon}(u)-2^*|u|_{2^*}^{2^*}\\
=&-\Big[\big[\sum_{i=1}^{l}\big(2^*(s_i)-2\big)b_{i,\varepsilon}(u)\big]+(2^*-2)|u|_{2^*}^{2^*}\Big]\nonumber\\
\leq&-\big[2^*(s_l)-2\big] \big[\sum_{i=1}^{l} b_{i,\varepsilon}(u)+|u|_{2^*}^{2^*}\big]\nonumber\\
=&-\big[2^*(s_l)-2\big] \|u\|^2
\leq-\big[2^*(s_l)-2\big] \delta_\varepsilon^2<0.\nonumber
\end{align}
 However, when $k\neq l$, we note that $b_{i,\varepsilon}(u)>0$ for $1\leq i\leq k$ and  $b_{i,\varepsilon}(u)>0$ for $k+1\leq i\leq l$. Here we view $b_{0,\varepsilon}(u)$ as $|u|_{2^*}^{2^*}$. Hence, we have
\begin{align}\lab{2014-12-24-we2}
\left\langle J'_\varepsilon(u), u\right\rangle =&2\|u\|^2-\sum_{i=0}^{l} 2^*(s_i)b_{i,\varepsilon}(u)\\
<&2\|u\|^2-2^*(s_{k})\sum_{i=0}^{l} b_{i,\varepsilon}(u)\nonumber\\
=&-\left[2^*(s_{k})-2\right]\|u\|^2
\leq -\left[2^*(s_{k})-2\right]\delta_\varepsilon^2<0.\nonumber
\end{align}
Apply the similar arguments as Lemma \ref{2014-12-16-l2}, we see that $\{u_n\}$ is bounded in $D_{0}^{1,2}(\R^N)$.
Let $\{t_n\}\subset \R$ be a sequence of multipliers satisfying
\be\lab{2014-12-16-e14}
\Phi'_{\varepsilon}(u_n)=\Phi'_{\varepsilon}\big|_{\mathcal{N}_\varepsilon}(u_n)+t_n J'_\varepsilon(u_n).
\ee
Testing by $u_n$, we obtain that
\be\lab{2014-12-16-e15}
t_n \langle J'_{\varepsilon}(u_n), u_n\rangle\rightarrow 0.
\ee
By \eqref{2014-12-16-e13} or \eqref{2014-12-24-we2} and \eqref{2014-12-16-e15}, we see that
\be\lab{2014-12-16-e16}
t_n\rightarrow 0\;\hbox{as}\;n\rightarrow +\infty.
\ee
Noting that $J'_{\varepsilon}(u_n)$ is bounded due to the boundedness of $\{u_n\}$, hence by \eqref{2014-12-16-e14} and \eqref{2014-12-16-e16}, we have  $
\Phi'_\varepsilon(u_n)\rightarrow 0\;\hbox{as}\;n\rightarrow +\infty. $
\ep
\br\lab{2014-12-16-r1}
By the formula \eqref{2014-12-16-e11}, for $\varepsilon>0$, we have that
\be\lab{2014-12-16-e18}
c_\varepsilon\geq \left[\frac{1}{2}-\frac{1}{2^*(s_k)}\right]\delta_\varepsilon^2>0.
\ee
Especially, when  $k=l$, by Lemma \ref{2014-12-16-l3}, we have
\be\lab{2014-12-24-we3}
c_\varepsilon>\left[\frac{1}{2}-\frac{1}{2^*(s_l)}\right]\delta_0^2>0.
\ee
For the case of $k\neq l$, we will prove that $c_\varepsilon$ is also achieved by some $u_\varepsilon$ and $u_\varepsilon$ is a mountain pass type solution (see Theorem \ref{2014-12-16-th1}). Set
\be\lab{2014-12-25-xe3}
\tilde{\Phi}_\varepsilon(u):=\frac{1}{2}\|u\|^2-\sum_{i=1}^{k}\int_{\R^N}\frac{\lambda_i}{2^*(s_i)}a_{i,\varepsilon}(x)|u|^{2^*(s_i)}.
\ee
It follows that there exists some $\tilde{\delta}_0>0$ such that
\be\lab{2014-12-25-xe4}
\inf_{u\in D_{0}^{1,2}(\R^N)\backslash \{0\}}\max_{t>0}\tilde{\Phi}_\varepsilon(tu)\geq \left[\frac{1}{2}-\frac{1}{2^*(s_{k})}\right]\tilde{\delta}_0^2>0.
\ee
Then it follows that
\begin{align}\lab{2014-12-25-xe5}
c_\varepsilon=&\Phi_{\varepsilon}(u_\varepsilon)=\max_{t>0}\Phi_{\varepsilon}(tu_\varepsilon)
>\max_{t>0}\tilde{\Phi}_{\varepsilon}(tu_\varepsilon)\geq\tilde{c}_\varepsilon
\geq \left[\frac{1}{2}-\frac{1}{2^*(s_{k})}\right]\tilde{\delta}_0^2>0.
\end{align}
\er
\bl\lab{2014-12-16-l5}
If  $k=l$,
$c_\varepsilon$ is strictly increasing by $\varepsilon$ in $[0,s_1)$.
\el
\bp
Let $0\neq u\in D_{0}^{1,2}(\R^N)$ be fixed.
For any $\varepsilon\in [0,s_1)$,
by Lemma \ref{2014-12-16-l1}, there exists $t_\varepsilon>0$ such that $t_{\varepsilon}u\in \mathcal{N}_{\varepsilon}$ and $t_\varepsilon$ is implicity given by
\be\lab{2014-12-16-e19}
\|u\|^2-\sum_{i=1}^{l} b_{i,\varepsilon}(u)t_\varepsilon^{2^*(s_i)-2}-|u|_{2^*}^{2^*}t_\varepsilon^{2^*-2}=0
\ee
By the Implicit Function Theorem, we see that $t(\varepsilon)\in C^1(\R)$ and $\frac{d}{d\varepsilon}t(\varepsilon)>0$ due to Lemma \ref{2014-12-16-l1}. Hence, recalling that $b_{i,\varepsilon}(u)$ is strictly decreasing by $\varepsilon$ and the formula \eqref{2014-12-16-e19}, we have
\begin{align}\lab{2014-12-16-e20}
&\frac{d}{d\varepsilon}\Phi_\varepsilon(t_\varepsilon u)\\
=&\frac{d}{d\varepsilon}\Big[\frac{1}{2}\|u\|^2 t_\varepsilon^2-\sum_{i=1}^{l}\frac{b_{i,\varepsilon}(u)}{2^*(s_i)}t_{\varepsilon}^{2^*(s_i)}-\frac{1}{2^*}|u|_{2^*}^{2^*}t_{\varepsilon}^{2^*}\Big]\nonumber\\
=&\frac{t'_\varepsilon}{t_\varepsilon}\big[\|u\|^2t_\varepsilon^2-\sum_{i=1}^{l}b_{i,\varepsilon}(u)t_{\varepsilon}^{2^*(s_i)}-|u|_{2^*}^{2^*}t_{\varepsilon}^{2^*}\big]-\sum_{i=1}^{l}\frac{b'_{i,\varepsilon}(u)}{2^*(s_i)}t_{\varepsilon}^{2^*(s_i)}\nonumber\\
=&-\sum_{i=1}^{l}\frac{b'_{i,\varepsilon}(u)}{2^*(s_i)}t_{\varepsilon}^{2^*(s_i)}>0,\nonumber
\end{align}
therefore,  $c_\varepsilon$ is strictly increasing by $\varepsilon$ in $[0,s_1)$.
\ep

\subsection{Existence of positive ground state of the approximating problem \eqref{Pva}}
In this subsection, we assume that $\varepsilon\in (0,s_1)$ is  fixed.

\bt\lab{2014-12-16-th1}
Let $N\geq 3,  0<s_1<s_2<\cdots<s_l<2$.  Suppose that there exists some $1\leq k\leq l$ such that $\lambda_i>0$ for $1\leq i\leq k$ and $\lambda_i<0$ for $k+1\leq i\leq l$. Furthermore, if $N=3$ and $k\neq l$, we assume that either  $s_1<1$ or $1\leq s_1<2$ with $\displaystyle\max\left\{|\lambda_{k+1}|,\cdots,|\lambda_l|\right\}$ small enough. Then for any  $\varepsilon\in (0,s_1)$, problem \eqref{Pva} possesses a positive ground state solution having the least energy
\be\lab{2014-12-16-e21}
c_\varepsilon<\frac{1}{N}S^{\frac{N}{2}}.
\ee
In particular, if $k=l$,    $c_\varepsilon$ is  increasing  strictly by $\varepsilon$.
\et

We postpone the proof of Theorem \ref{2014-12-16-th1} and do a little preparation before that. Denote
\be\lab{2014-12-16-e22}
\Psi(u):=\frac{1}{2}\|u\|^2-\frac{1}{2^*}|u|_{2^*}^{2^*}\;\;u\in D_{0}^{1,2}(\R^N).
\ee
It is well known that
\be\lab{2014-12-16-e23}
\min_{u\in D_{0}^{1,2}(\R^N)\backslash\{0\}}\max_{t>0}\Psi(tu)=\frac{1}{N}S^{\frac{N}{2}}.
\ee
By Lemma \ref{2014-12-16-l2}, any $(PS)_c$ sequence of $\Phi_\varepsilon$ is bounded in $D_{0}^{1,2}(\R^N)$. Hence, we may give the following proposition:
\bo\lab{2014-12-16-prop1}
Let $N\geq 3, 0<s_1<s_2<\cdots<s_l<2$, there exists some $1\leq k\leq l$ such that $\lambda_i>0$ for $1\leq i\leq k$ and $\lambda_i<0$ for $k+1\leq i\leq l$. Take $\varepsilon\in (0,s_1)$ and assume that $\{u_n\}$ is a $(PS)_c$ sequence of $\Phi_\varepsilon$, i.e.,
\be\lab{2014-12-16-e24}
\begin{cases}
\Phi_\varepsilon(u_n)\rightarrow c,\\
\Phi'_\varepsilon(u_n)\rightarrow 0.
\end{cases}
\ee
Up to a subsequence, we may assume that $u_n\rightharpoonup u_0$ in $D_{0}^{1,2}(\R^N)$ and $u_n\rightarrow u_0$ a.e. in $\R^N$. Denote $\tilde{u}_n:=u_n-u_0$, then $\{\tilde{u}_n\}$ is a $PS$ sequence of $\Psi$.
\eo
\bp
Basing on the Lemma \ref{2014-12-12-wl1}, by H\"older inequality and Hardy Sobolev inequality, it is easy to prove that
\be\lab{2014-12-16-e25}
\int_{\R^N}a_{i,\varepsilon}(x)|u_n|^{2^*(s_i)-1}hdx -\int_{\R^N}a_{i,\varepsilon}(x)|u_0|^{2^*(s_i)-1}hdx=o(1)\|h\|.
\ee
Since $\{u_n\}$ is a $(PS)_c$ sequence of $\Phi_\varepsilon$, we see that $\Phi'_\varepsilon(u_0)=0$. Then it follows that
\be\lab{2014-12-16-e26}
\int_{\R^N}\nabla (u_n-u_0)\nabla hdx-\int_{\R^N}\big[|u_n|^{2^*-2}u_n -|u_0|^{2^*-2}u_0\big]hdx=o(1)\|h\|.
\ee
By the Br\'ezis-Lieb Lemma, we see that
\be\lab{2014-12-16-e27}
|u_n|^{2^*-2}u_n -|u_0|^{2^*-2}u_0 -|u_n-u_0|^{2^*-2}(u_n-u_0) \rightarrow 0\;\hbox{strongly in}\;L^{\frac{2^*}{2^*-1}}(\R^N).
\ee
Hence, by \eqref{2014-12-16-e26} and \eqref{2014-12-16-e27}, this proposition is proved.
\ep

\bc\lab{2014-12-16-cro1}
Let $N\geq 3, 0<s_1<s_2<\cdots<s_l<2$, there exists some $1\leq k\leq l$ such that $\lambda_i>0$ for $1\leq i\leq k$ and $\lambda_i<0$ for $k+1\leq i\leq l$. Then for any $\varepsilon\in (0,s_1)$, $\Phi_\varepsilon$ satisfies $(PS)_c$ condition
if $c<\frac{1}{N}S^{\frac{N}{2}}$.
\ec
\bp
Let $\{u_n\}$ be a $(PS)_c$ sequence of $\Phi_\varepsilon$ with $c<\frac{1}{N}S^{\frac{N}{2}}$. Up to a subsequence, we assume that $u_n\rightharpoonup u_0$ in $D_{0}^{1,2}(\R^N)$ and $u_n\rightarrow u_0$ a.e. in $\R^N$. We prove this corollary by the way of negation. If $u_n\not\rightarrow u_0$ in $D_{0}^{1,2}(\R^N)$, then by Proposition \ref{2014-12-16-prop1}, we see that $\tilde{u_n}:=u_n-u_0$ is a $PS$ sequence of $\Psi$ and $\tilde{u}_n\not\rightarrow 0$. Then it is easy to prove that
\be\lab{2014-12-16-e28}
\lim_{n\rightarrow +\infty}\Psi(\tilde{u}_n)\geq \frac{1}{N}S^{\frac{N}{2}}.
\ee
By Br\'ezis-Lieb Lemma again, we have that
\be\lab{2014-12-16-e29}
\Phi_\varepsilon(u_n)=\Phi_\varepsilon(u_0)+\Psi(\tilde{u}_n)+o(1).
\ee
Since $\Phi'_\varepsilon(u_0)=0$, it is easy to see that $\Phi_\varepsilon(u_0)\geq 0$. Hence, by \eqref{2014-12-16-e28} and \eqref{2014-12-16-e29}, we obtain that
\be\lab{2014-12-16-e30}
c=\lim_{n\rightarrow +\infty}\Phi_\varepsilon(u_n)\geq \frac{1}{N}S^{\frac{N}{2}},
\ee
a contradiction.
Hence, $u_n\rightarrow u_0$ strongly in $D_{0}^{1,2}(\R^N)$ and it follows that $\Phi_\varepsilon(u_0)= c$.
\ep

\bl\lab{2014-12-16-l6}
Let $N\geq 3, 0<s_1<s_2<\cdots<s_l<2$. Suppose that  there exists some $1\leq k\leq l$ such that $\lambda_i>0$ for $1\leq i\leq k$ and $\lambda_i<0$ for $k+1\leq i\leq l$. Furthermore, if $N=3$ and $k\neq l$, we assume that either  $s_1<1$ or $1\leq s_1<2$ with $\displaystyle\max\left\{|\lambda_{k+1}|,\cdots,|\lambda_l|\right\}$ small enough. Take $\varepsilon\in [0,s_1)$ and suppose that $c_\varepsilon$ is given by \eqref{2014-12-5-e3}, then we have
\be\lab{2014-12-16-e31}
c_\varepsilon<\frac{1}{N}S^{\frac{N}{2}}.
\ee
\el
\bp
Let $\displaystyle U(x):=\frac{[N(N-2)]^{\frac{N-2}{4}}}{[1+|x|^2]^{\frac{N-2}{2}}}$.
For the case of $k=l$, it is easy to see that
\begin{align}\lab{2014-12-26-ze1}
c_\varepsilon\leq \max_{t>0}\Phi_{\varepsilon}(tU)<\max_{t>0}\Psi(tU)=\frac{1}{N}S^{\frac{N}{2}}.
\end{align}
Next, we consider the case of $k\neq l$. When $1\leq s_1<2$ with $\displaystyle\max\left\{|\lambda_{k+1}|,\cdots,|\lambda_l|\right\}$ small enough, a direct computation shows  that \eqref{2014-12-26-ze1} is also satisfied. And we note that the small bound can be chosen independent of $\varepsilon$ for $\varepsilon$ small enough.
When $0<s_1<1$ and $k\neq l$, we let $0\neq x_0\in \R^N$, and $\psi\in D_{0}^{1,2}(\R^N)$ be a nonnegative function such that $\psi\equiv 1$ on $B(0,\rho), 0<\rho<|x_0|$.
For $\sigma >0$, we define
$$U_\sigma(x):=\sigma^{\frac{2-N}{2}}U\left(\frac{x-x_0}{\sigma}\right),u_\sigma(x):=\psi(x-x_0)U_\sigma(x).$$
Noting that $\varepsilon>0$ is fixed, a direct computation shows that
\be\lab{2014-12-24-we1}
\int_{\R^N}a_{i,\varepsilon}(x)|u_\sigma(x)|^{2^*(s_i)}dx=O\left(\sigma^{s_i}\right),  \;i=1,2,\cdots,l,
\ee
and
\be\lab{2014-12-25-ze1}
\int_{\R^N}|\nabla u_\sigma(x)|^2dx=S^{\frac{N}{2}}+O(\sigma^{N-2});\int_{\R^N}|u_\sigma(x)|^{2^*}dx=S^{\frac{N}{2}}+O(\sigma^N).
\ee
 Since $\lambda_1>0$ and $s_1<\min\{s_2,\cdots, s_l, N-2\}$ under the assumptions, then it is easy to see that

\be\lab{2014-12-16-e32}
\max_{t>0}\Phi_{\varepsilon}(tu_\sigma)<\frac{1}{N}S^{\frac{N}{2}}\;\hbox{for $\sigma$ small enough}.
\ee
Hence, by the definition of $c_\varepsilon$, we obtain that
$\displaystyle c_\varepsilon<\frac{1}{N}S^{\frac{N}{2}}$.
\ep

\vskip 0.2in
\noindent{\bf Proof of Theorem \ref{2014-12-16-th1}:} Let $\{u_n\}\subset \mathcal{N}_\varepsilon$ be a minimizing sequence of $c_\varepsilon$. Then by Lemma \ref{2014-12-16-l4}, we see that $\{u_n\}$ is also a $(PS)_{c_\varepsilon}$ sequence of $\Phi_\varepsilon$. Under the assumptions, by Lemma \ref{2014-12-16-l6}, we have $\displaystyle c_\varepsilon<\frac{1}{N}S^{\frac{N}{2}}$. By Corollary \ref{2014-12-16-cro1}, we observe  that $\Phi_\varepsilon$ satisfies the $(PS)_{c_\varepsilon}$ condition. Up to a subsequence, we may assume that $u_n\rightarrow u_0$ strongly in $D_{0}^{1,2}(\R^N)$ and $\Phi_\varepsilon(u_0)=c_\varepsilon$. Hence, $u_0$ is a minimizer of $c_\varepsilon$. Noting that $\Phi_\varepsilon$ is even, we see that $|u_0|$ is also a minimizer of $c_\varepsilon$. Hence, without loss of generality, we may assume that $u_0\geq 0$. Then, we see that $\Phi'_\varepsilon(u_0)=0$. By the maximum principle, we have that $u_0>0$ in $\R^N\backslash\{0\}$. Hence, $u_0$ is a positive ground state solution of problem \eqref{Pva} and Theorem \ref{2014-12-16-th1} is proved.  \hfill $\Box$

\vskip 0.2in
\br\lab{2014-12-16-xr1}
For $\varepsilon\in [0,s_1)$,
we define the mountain pass  value
\be\lab{2014-3-10-e2}
\tilde{c}_\varepsilon:=\inf_{\gamma\in \Gamma_\varepsilon}\max_{t\in [0,1]}\Phi_{\varepsilon}(\gamma(t)),
\ee
 where $\displaystyle \Gamma_\varepsilon:=\left\{\gamma(t)\in C\left([0,1],  D_{0}^{1,2}(\R^N)\right): \gamma(0)=0, \Phi_{\varepsilon}(\gamma(1))<0\right\}$.
It is standard to prove that $c_\varepsilon=\tilde{c}_\varepsilon$ and any ground state solution of $\eqref{Pva}$ is a mountain pass solution provided that $\varepsilon>0$. Precisely, by the definition and the result of Lemma \ref{2014-12-16-l1}, it is easy to see that $\tilde{c}_\varepsilon\leq c_\varepsilon$ for all $\varepsilon\in [0,s_1)$. When $\varepsilon>0$, by Corollary \ref{2014-12-16-cro1}, $\Phi_\varepsilon$ satisfies $(PS)_{\tilde{c}_\varepsilon}$ condition. Hence, there exists a mountain pass solution $\tilde{u}_\varepsilon$  such that  $\Phi_\varepsilon(\tilde{u}_\varepsilon)=\tilde{c}_\varepsilon$.  It follows that
  \be
  \tilde{c}_\varepsilon=\Phi_\varepsilon(\tilde{u}_\varepsilon)
  =\max_{t>0}\Phi_\varepsilon(t\tilde{u}_\varepsilon)\geq \min_{u\neq 0}\max_{t>0}\Phi_\varepsilon(tu)=c_\varepsilon,
  \ee
  thus we can obtain the reverse inequality. Hence, we have $\tilde{c}_\varepsilon=c_\varepsilon$ for $\varepsilon\in(0,s_1)$.
\er

\bl\lab{2014-12-16-l7}
Under the assumptions of Theorem \ref{2014-12-12-mainth2}, we have that
$\displaystyle \limsup_{\varepsilon\rightarrow 0^+}\tilde{c}_\varepsilon\leq \tilde{c}_0$. Moreover, if $k=l$, we have that
$\tilde{c}_{\varepsilon}>\tilde{c}_0$ and thus $\displaystyle\lim_{\varepsilon\rightarrow 0^+} \tilde{c}_{\varepsilon}=\tilde{c}_0$.
\el
\bp
For any $\delta>0$, there exists $\gamma_0\in \Gamma_0$ such that
\be\lab{2014-12-16-xe1}
\max_{t\in [0,1]}\Phi_0(\gamma_0(t))<\tilde{c}_0+\delta.
\ee
Denote $\gamma_0(1)=\phi$, since $\gamma_0\in \Gamma_0$,  we have $\Phi_0(\phi)<0$. By the Lebesgue's dominated convergence theorem, we have
\be\lab{2014-12-16-xe2}
\lim_{\varepsilon\rightarrow 0}\int_{\R^N}\lambda_i a_{i,\varepsilon}(x)|\phi|^{2^*(s_i)}dx=\int_{\R^N}\lambda_i \frac{1}{|x|^{s_i}}|\phi|^{2^*(s_i)}dx.
\ee
Hence, by the continuity, we see that $\gamma_0\in \Gamma_\varepsilon$ when $\varepsilon$ is small enough. Now, take $\varepsilon_n \downarrow 0$ and denote $t_n\in [0,1]$ such that
\be\lab{2014-12-16-xe3}
\Phi_{\varepsilon_n}(\gamma_0(t_n))=\max_{t\in [0,1]}\Phi_{\varepsilon_n}(\gamma_0(t)).
\ee
Up to a subsequence if necessary, we may assume that $t_n\rightarrow t^*\in [0,1]$. Set $u_n:=\gamma_0(t_n)$ and $u^*:=\gamma_0(t^*)$, since $\gamma_0\in C([0,1], D_{0}^{1,2}(\R^N))$, we obtain that $u_n\rightarrow u^*$ strongly in $D_{0}^{1,2}(\R^N)$. Hence, we have
\be\lab{2014-12-16-xe4}
\Phi_{\varepsilon_n}(u_n)=\Phi_{\varepsilon_n}(u^*)+o(1).
\ee
On the other hand, by the Lebesgue's dominated convergence theorem again, we have
\be\lab{2014-12-16-xe5}
\Phi_{\varepsilon_n}(u^*)=\Phi_{0}(u^*)+o(1).
\ee
Then by \eqref{2014-12-16-xe4} and \eqref{2014-12-16-xe5}, we have
\begin{align}\lab{2014-12-16-xe6}
\tilde{c}_\varepsilon\leq&\Phi_{\varepsilon_n}(\gamma_0(t_n))=\Phi_{\varepsilon_n}(u_n)\nonumber\\
=&\Phi_0(u^*)+o(1)\leq \max_{t\in[0,1]}\Phi_0(\gamma_0(t))+o(1)\nonumber\\
\leq&\tilde{c}_0+\delta+o(1).
\end{align}
Hence, $\displaystyle \limsup_{n\rightarrow +\infty}\tilde{c}_{\varepsilon_n}\leq \tilde{c}_0$ due to the arbitrariness  of $\delta$.

\vskip 0.02in
Moreover, if $k=l$, noting that
$a_{i,\varepsilon}(x)$ is decreasing by $\varepsilon$ for all $i=1,2,\cdots,l$, it is easy to see that $\tilde{c}_\varepsilon\geq \tilde{c}_0$. By Theorem \ref{2014-12-16-th1} and Remark \ref{2014-12-16-xr1}, $\tilde{c}_\varepsilon=c_\varepsilon$ can be achieved. Hence, $\tilde{c}_\varepsilon>\tilde{c}_0$.
It is also trivial that  $\displaystyle \lim_{\varepsilon\rightarrow 0^+} \tilde{c}_{\varepsilon}\geq \tilde{c}_0$. Hence,
\be\lab{2014-12-16-xe7}
\lim_{\varepsilon\rightarrow 0^+} \tilde{c}_{\varepsilon}= \tilde{c}_0.
\ee
\ep

\bl\lab{2014-12-16-l8}
$\tilde{c}_0\leq c_0$ and thus $\displaystyle \limsup_{\varepsilon\rightarrow 0^+} c_{\varepsilon}\leq c_0$. Especially, $\displaystyle \tilde{c}_0=c_0\;\hbox{and}\;\lim_{\varepsilon\rightarrow 0^+}c_{\varepsilon}=c_0$ provided $k=l$.
\el
\bp
For any $0\neq u\in D_{0}^{1,2}(\R^N)$,
since  $\gamma(t):=tTu\in \Gamma_0(t)$ for $T$ large enough, then by the definition of $c_0$ and $\tilde{c}_0$, it is easy to see that
\be\lab{2014-12-16-xe10}
\tilde{c}_0\leq c_0.
\ee
On the other hand, by Remark \ref{2014-12-16-xr1} and Lemma \ref{2014-12-16-l7}, we have
\be\lab{2014-12-16-xe8}
\tilde{c}_0\geq \limsup_{\varepsilon\rightarrow 0^+}\tilde{c}_\varepsilon=\limsup_{\varepsilon\rightarrow 0^+}c_\varepsilon.
\ee
Moreover, if $k=l$,
by Lemma \ref{2014-12-16-l5}, $c_\varepsilon>c_0$ for any $\varepsilon>0$, combining with Lemma \ref{2014-12-16-l7}, we obtain the reverse inequality
\be\lab{2014-12-16-xe9}
\tilde{c}_0=\lim_{\varepsilon\rightarrow 0^+}\tilde{c}_\varepsilon=\lim_{\varepsilon\rightarrow 0^+}c_\varepsilon\geq c_0.
\ee
Hence, by \eqref{2014-12-16-xe10} and \eqref{2014-12-16-xe9}, we see that $c_0=\tilde{c}_0$.
\ep
\br\lab{2014-12-25-zr1}
When $\varepsilon=0$, since it is not trivial to see that $c_0$ is a ground state value, we can not obtain that $\tilde{c}_0=c_0$ by the arguments as the case of $\varepsilon>0$ that mentioned in Remark \ref{2014-12-16-xr1}. However, if $k=l$, by Lemma \ref{2014-12-16-l8} above, we still obtain that $\tilde{c}_0=c_0$. For the case of $k\neq l$, since we can not obtain the monotonicity of $c_\varepsilon$, we are unable to  get the conclusion of $\tilde{c}_0=c_0$ up to now. However, we note that after the results established in present paper, we will see that this relationship still holds. Especially,  $c_0$ can be attained. We can also obtain that $\displaystyle \lim_{\varepsilon\rightarrow 0^+}c_{\varepsilon}=c_0$ for the case $k\neq l$.
\er

\s{Interpolation Inequalities and Pohozaev Identity}
\renewcommand{\theequation}{4.\arabic{equation}}
\renewcommand{\theremark}{4.\arabic{remark}}
\renewcommand{\thedefinition}{4.\arabic{definition}}
The following Propositions \ref{2014-5-5-interpolation-corollary}-\ref{2014-4-22-interpolation-corollary} are proved in  \cite{ZhongZou.2017} and Proposition \ref{2014-6-21-prop1}  is obtained in  \cite{LehrerMaia.2013}.    Define
\be\lab{2014-5-8-we1}
\vartheta(s_1,s_2):=\frac{N(s_2-s_1)}{s_2(N-s_1)}\;\quad \hbox{for \;$0\leq s_1\leq s_2\leq 2. $}
\ee
\bo\lab{2014-5-5-interpolation-corollary} (see \cite[Corollary 1]{ZhongZou.2017})  Let $\Omega\subset \R^N (N\geq 3)$ be an open set.
Assume  $0\leq s_1<2$.  Then for any $s_2\in [s_1,2]$ and $\theta\in [\vartheta(s_1,s_2), 1]$, there exists $C(\theta)>0$ such that
\be\lab{2014-5-8-e1}
|u|_{2^*(s_1), {s_1}}\leq C(\theta)\|u\|^\theta |u|_{2^*(s_2),  {s_2}}^{1-\theta}
\ee
for all $u\in D_{0}^{1,2}(\Omega)$, where $\|u\|:=\big(\int_\Omega |\nabla u|^2dx\big)^{\frac{1}{2}}$.
\eo

Define
\be\lab{2014-5-8-we2}
\varsigma(s_1,s_2):=\frac{(N-s_1)(2-s_2)}{(N-s_2)(2-s_1)}\;\quad \hbox{for   \;$0\leq s_1\leq s_2\leq 2.$}
\ee
\bo\lab{2014-4-22-interpolation-corollary} (see \cite[Corollary 2]{ZhongZou.2017})  Let $\Omega\subset \R^N (N\geq 3)$ be an open set.
Assume  $0<s_2\leq 2$.   Then for any $s_1\in [0,s_2]$ and $\sigma\in [0,\varsigma(s_1,s_2)]$, there exists a  $C(\sigma)>0$ such that
\be\lab{2014-5-8-we3}
|u|_{2^*(s_2),   {s_2}}\leq C(\sigma)\|u\|^{1-\sigma} |u|_{2^*(s_1),  {s_1}}^{\sigma}
\ee
for all $u\in D_{0}^{1,2}(\Omega)$
\eo

\bo\lab{2014-6-21-prop1}(see \cite[Proposition 2.1]{LehrerMaia.2013})
Let $u\in H^1(\Omega)\backslash \{0\}$ be a solution to  the  equation $-\Delta u=g(x, u)$ and $G(x, u)=\int_0^u g(x, s)ds$ is such that $G\big(\cdot, u(\cdot)\big)$ and $x_iG_{x_i}\big(\cdot, u(\cdot)\big)$ are in $L^1(\Omega)$, then $u$ satisfies:
$$\int_{\partial\Omega}|\nabla u|^2 x\cdot \eta dS_x=2N\int_\Omega G(x, u)dx+2\sum_{i=1}^{N}\int_\Omega x_iG_{x_i}(x, u)dx-(N-2)\int_\Omega|\nabla u|^2dx,$$
where $\Omega$ is a regular domain in $\R^N$ and $\eta$ denotes the unitary exterior normal vector to $\partial\Omega$. Moreover, if $\Omega=\R^N$, then
$$2N\int_{\R^N}G(x, u)dx+2\sum_{i=1}^{N}\int_{\R^N}x_iG_{x_i}(x, u)dx=(N-2)\int_{\R^N}|\nabla u|^2 dx.$$
\eo

\bc\lab{2014-12-25-cro1}
For $\varepsilon>0$ small enough, we still have that for any $\theta\in [\vartheta(s_i,s_j), 1]$ if $0<s_i\leq s_j< 2$,
\be\lab{2014-12-25-xe1}
\left(\int_{\R^N} a_{i,\varepsilon}(x)|u|^{2^*(s_i)}dx\right)^{\frac{1}{2^*(s_i)}}\leq C(\theta)\|u\|^\theta \left(\int_{\R^N} a_{j,\varepsilon}(x)|u|^{2^*(s_j)}dx\right)^{\frac{1-\theta}{2^*(s_j)}}
\ee
And for any $\sigma\in [0,\varsigma(s_i,s_j)]$ if $0< s_i\leq s_j<2$,
\be\lab{2014-12-25-xe2}
\left(\int_{\R^N} a_{j,\varepsilon}(x)|u|^{2^*(s_j)}dx\right)^{\frac{1}{2^*(s_j)}}\leq C(\sigma)\|u\|^{1-\sigma} \left(\int_{\R^N} a_{i,\varepsilon}(x)|u|^{2^*(s_i)}dx\right)^{\frac{\sigma}{2^*(s_i)}}.
\ee
\ec
\bp
We replace $dx$ by the  new measure $d\nu:=\begin{cases} \frac{dx}{|x|^{-\varepsilon}}\;\hbox{if}\;|x|\leq 1\\  \frac{dx}{|x|^{\varepsilon}}\;\hbox{if}\;|x|> 1\end{cases}$. Recalling the embedding relationship in Lemma \ref{2014-12-12-wl1}, by the same arguments as  the the proofs of \cite[Corollary 1 and Corollary 2]{ZhongZou.2017}, we can obtain the results of \eqref{2014-12-25-xe1} and \eqref{2014-12-25-xe2}. We omit the details.
\ep

\bc\lab{2014-12-16-cro2}
Let $N\geq 3, 0<s_i<2$ and $\varepsilon\in (0,s_1)$. Then any solution of \eqref{Pva} satisfies
\be\lab{2014-12-16-e33}
\int_{\mathbb{B}_1} \sum_{i=1}^{l}\frac{\lambda_i a_{i,\varepsilon}(x)}{2^*(s_i)}|u|^{2^*(s_i)}dx=\int_{\mathbb{B}_1^c} \sum_{i=1}^{l}\frac{\lambda_i a_{i,\varepsilon}(x)}{2^*(s_i)}|u|^{2^*(s_i)}dx.
\ee
\ec
\bp
Take $\displaystyle G(x, u)=\sum_{i=1}^{l}\frac{1}{2^*(s_i)}\lambda_i a_{i,\varepsilon}(x)|u|^{2^*(s_i)}+\frac{1}{2^*}u^{2^*}$.
By Proposition  \ref{2014-6-21-prop1}, we have
\begin{align}\lab{2014-12-17-e1}
&2N\int_{\R^N}\Big[\sum_{i=1}^{l}\frac{\lambda_i}{2^*(s_i)}a_{i,\varepsilon}(x)|u|^{2^*(s_i)}+\frac{1}{2^*}|u|^{2^*}\Big]dx\\
&+2\sum_{j=1}^{N}\int_{\R^N}\sum_{i=1}^{l}\frac{\lambda_i}{2^*(s_i)}\frac{\partial}{\partial x_j}a_{i,\varepsilon}(x)|u|^{2^*(s_i)}x_j\nonumber\\
=&(N-2)\int_{\R^N}|\nabla u|^2dx.\nonumber
\end{align}
Noting that
\be\lab{2014-12-17-e2}
\frac{\partial}{\partial x_j}a_{i,\varepsilon}(x)=\begin{cases}
-(s_i-\varepsilon)\frac{1}{|x|^{s_i+2-\varepsilon}}x_j\quad &\hbox{for}\;|x|<1,\\
-(s_i+\varepsilon)\frac{1}{|x|^{s_i+2+\varepsilon}}x_j\quad&\hbox{for}\;|x|>1,
\end{cases}
\ee
we obtain that
\be\lab{2014-12-17-e3}
\sum_{j=1}^{N}\frac{\partial}{\partial x_j}a_{i,\varepsilon}(x)x_j=\begin{cases}
-(s_i-\varepsilon)a_{i,\varepsilon}(x),\quad &|x|<1,\\
-(s_i+\varepsilon)a_{i,\varepsilon}(x), \quad &|x|>1.
\end{cases}
\ee
Then, substitute into \eqref{2014-12-17-e1}, we obtain that
\begin{align}\lab{2014-12-17-e4}
&2N\int_{\R^N}\Big[\sum_{i=1}^{l}\frac{\lambda_i}{2^*(s_i)}a_{i,\varepsilon}(x)|u|^{2^*(s_i)}+\frac{1}{2^*}|u|^{2^*}\Big]dx
-\int_{\R^N}\Big[\sum_{i=1}^{l}\frac{2s_i\lambda_i}{2^*(s_i)}a_{i,\varepsilon}(x)|u|^{2^*(s_i)}\Big]dx\\
&+2\varepsilon \int_{\mathbb{B}_1}\Big[\sum_{i=1}^{l}\frac{\lambda_ia_{i,\varepsilon}(x)}{2^*(s_i)}|u|^{2^*(s_i)}\Big]dx
-2\varepsilon \int_{\mathbb{B}_1^c}\Big[\sum_{i=1}^{l}\frac{\lambda_ia_{i,\varepsilon}(x)}{2^*(s_i)}|u|^{2^*(s_i)}\Big]dx\nonumber\\
=&(N-2)\int_{\R^N}|\nabla u|^2 dx.\nonumber
\end{align}
On the other hand, since $u$ is a solution, we have
\be\lab{2014-12-17-e5}
\int_{\R^N}|\nabla u|^2 dx=\int_{\R^N}\Big[\sum_{i=1}^{l}\lambda_i a_{i,\varepsilon}(x)|u|^{2^*(s_i)}+|u|^{2^*}\Big]dx.
\ee
Hence, by \eqref{2014-12-17-e4} and \eqref{2014-12-17-e5}, we get
\be\lab{2014-12-17-e6}
\int_{\mathbb{B}_1}\Big[\sum_{i=1}^{l}\frac{\lambda_ia_{i,\varepsilon}(x)}{2^*(s_i)}|u|^{2^*(s_i)}\Big]dx
= \int_{\mathbb{B}_1^c}\Big[\sum_{i=1}^{l}\frac{\lambda_ia_{i,\varepsilon}(x)}{2^*(s_i)}|u|^{2^*(s_i)}\Big]dx.
\ee
\ep

\s{Proof of Theorem \ref{2014-12-12-mainth2}}
\renewcommand{\theequation}{5.\arabic{equation}}
\renewcommand{\theremark}{5.\arabic{remark}}
\renewcommand{\thedefinition}{5.\arabic{definition}}
\subsection{Preliminary}
\br\lab{2015-2-16-r1}
For $\forall\;\varepsilon\in (0,s_1)$, by Theorem \ref{2014-12-16-th1}, problem \eqref{Pva} possesses a positive ground state solution $u_\varepsilon$ such that $\Phi_\varepsilon(u_\varepsilon)=c_\varepsilon$. Now, we take $\varepsilon_n\downarrow 0$ as $n\rightarrow +\infty$ and assume that $u_n$ is a positive ground state solution of \eqref{Pva} with $\varepsilon=\varepsilon_n$.
Similar to the formula \eqref{2014-12-16-e11}, it is easy to prove that
\be\lab{2014-12-17-e7}
c_{\varepsilon_n}=\Phi_{\varepsilon_n}(u_n)\geq (\frac{1}{2}-\frac{1}{2^*(s_k)})\|u_n\|^2.
\ee
By Lemma \ref{2014-12-16-l8}, we see that $\displaystyle\limsup _{n\rightarrow +\infty}c_{\varepsilon_n}\leq c_0$. Hence, $\{u_n\}$ is bounded in $D_{0}^{1,2}(\R^N)$. Up to a subsequence, we assume that $u_n\rightharpoonup u_0$ in $D_{0}^{1,2}(\R^N)$ and $u_n\rightarrow u_0$ a.e. in $\R^N$.
\er

\bl\lab{2014-12-17-l1}
$u_0$ is a critical point of $\Phi_0$, i.e., $\Phi'_0(u_0)=0$.
\el
\bp
We claim that for any $\phi\in D_{0}^{1,2}(\R^N)$ and $i\in \{1,2,\cdots, l\}$, we have
\be\lab{2014-12-17-we1}
\lim_{n\rightarrow+\infty}\int_{\R^N} \left[ a_{i,\varepsilon_n}(x)|u_n|^{2^*(s_i)-2}u_n\phi\right]dx
=\int_{\R^N}\left[\frac{1}{|x|^{s_i}}|u_0|^{2^*(s_i)-2}u_0\phi\right]dx.
\ee
Without loss of generality, we may also assume that $\phi\geq 0$. Otherwise,  we write  $\phi=\phi_+-\phi_-$ and  discuss on $\phi_+$ and $\phi_-$,  respectively.
Firstly by the Fatou's Lemma, it is easy to see that
\be\lab{2014-12-17-e8}
\int_{\R^N}\left[\frac{1}{|x|^{s_i}}|u_0|^{2^*(s_i)-2}u_0\phi\right]dx\leq \liminf_{n\rightarrow +\infty}\int_{\R^N}\left[ a_{i,\varepsilon_n}(x)|u_n|^{2^*(s_i)-2}u_n\phi\right]dx.
\ee
On the other hand, since $a_{i,\varepsilon_n}(x)\leq a_{i,0}(x)=\frac{1}{|x|^{s_i}}$, we have
\be\lab{2014-12-17-we2}
\int_{\R^N}\left[ a_{i,\varepsilon_n}(x)|u_n|^{2^*(s_i)-2}u_n\phi\right]dx\leq \int_{\R^N}\left[ \frac{1}{|x|^{s_i}}|u_n|^{2^*(s_i)-2}u_n\phi\right]dx.
\ee
Since $u_n\rightharpoonup u_0$ in $D_{0}^{1,2}(\R^N)$, we see that
$$|u_n|^{2^*(s_i)-2}u_n\rightharpoonup |u_0|^{2^*(s_i)-2}u_0  \hbox{  in } L^{\frac{2^*(s_i)}{2^*(s_i)-1}}(\R^N).$$
Hence, we have
\be\lab{2014-12-17-we3}
\lim_{n\rightarrow +\infty}\int_{\R^N}\left[ \frac{1}{|x|^{s_i}}|u_n|^{2^*(s_i)-2}u_n\phi\right]dx=\int_{\R^N}\left[ \frac{1}{|x|^{s_i}}|u_0|^{2^*(s_i)-2}u_0\phi\right]dx.
\ee
By \eqref{2014-12-17-we2} and \eqref{2014-12-17-we3}, we obtain
\be\lab{2014-12-17-we4}
\limsup_{n\rightarrow +\infty}\int_{\R^N}\left[ a_{i,\varepsilon_n}(x)|u_n|^{2^*(s_i)-2}u_n\phi\right]dx\leq\int_{\R^N}\left[ \frac{1}{|x|^{s_i}}|u_0|^{2^*(s_i)-2}u_0\phi\right]dx.
\ee
Thus, \eqref{2014-12-17-we1} is proved by \eqref{2014-12-17-e8} and \eqref{2014-12-17-we4}.
Recalling  that  $u_n$ is a critical point of $\Phi_{\varepsilon_n}$, we have that
$
\langle \Phi'_{\varepsilon_n}(u_n), \phi\rangle =0\;\hbox{for all}\;\phi\in D_{0}^{1,2}(\R^N).
$
Then by \eqref{2014-12-17-we1} and $u_n\rightharpoonup u_0$ in $D_{0}^{1,2}(\R^N)$, we see that
$
\langle \Phi'_{0}(u_0), \phi\rangle =0\;\hbox{for all}\;\phi\in D_{0}^{1,2}(\R^N),
$ i.e.,  $\Phi'_0(u_0)=0$.
\ep

\bl\lab{2014-12-17-l2}
$\displaystyle 0\leq \Phi_0(u_0)\leq \lim_{n\rightarrow +\infty}c_{\varepsilon_n}\leq c_0$ and if $u_0\neq 0$, we have $\Phi_0(u_0)=c_0>0$.
\el
\bp
By Lemma \ref{2014-12-17-l1}, $\Phi'_0(u_0)=0$.
If $u_0=0$, we have $\Phi_0(u_0)=0$. If $u_0\neq 0$, it is easy to see that $u_0\in \mathcal{N}_0$, then it follows that $\Phi_0(u_0)\geq c_0>0$. Hence, we always have
\be\lab{2014-12-17-we7}
\Phi_0(u_0)\geq 0.
\ee
Since $u_n$ is a ground state solution of \eqref{Pva} with $\varepsilon=\varepsilon_n$, similar to \eqref{2014-12-16-e11}, we have that
\begin{align}\lab{2014-12-17-we8}
c_{\varepsilon_n}=\Phi_{\varepsilon_n}(u_n)
=&\left[\frac{1}{2}-\frac{1}{2^*(s_k)}\right]\|u_n\|^2
+\sum_{i=1}^{l}\left[\frac{1}{2^*(s_k)}-\frac{1}{2^*(s_i)}\right]\lambda_i |u_n|_{2^*(s_i),i,\varepsilon_n}^{2^*(s_i)} \nonumber\\ &+\left[\frac{1}{2^*(s_k)}-\frac{1}{2^*}\right]|u_n|_{2^*}^{2^*}.
\end{align}
Noting that $\displaystyle \left[\frac{1}{2^*(s_k)}-\frac{1}{2^*(s_i)}\right]\lambda_i>0$ for $i\neq k$, then by Fatou's Lemma and Lemma \ref{2014-12-16-l8}, we get that
\be\lab{2014-12-17-we9}
\Phi_0(u_0)\leq \liminf_{n\rightarrow +\infty}c_{\varepsilon_n}\leq c_0.
\ee
Furthermore, if $u_0\neq 0$, then by the definition of $c_0$, it is trivial to obtain the reverse inequality $\Phi_0(u_0)\geq c_0.$ Hence, $\Phi_0(u_0)=c_0$. Evidently, $c_0>0$, see also Remark \ref{2014-12-16-r1} and Lemma \ref{2014-12-16-l1}.
\ep

\bl\lab{2014-12-17-l3}
If $\displaystyle \lim_{n\rightarrow +\infty}\int_{\R^N}a_{i,\varepsilon_n}(x)|u_n|^{2^*(s_i)}dx=0$ for some $i\in \{1,2,\cdots,l\}$,
then $\{u_n\}$ is a $PS$ sequence of $\Psi$, i.e., $\Psi'(u_n)\rightarrow 0$.
\el
\bp
Noting that $\{u_n\}$ is bounded in $D_{0}^{1,2}(\R^N)$. By Corollary \ref{2014-12-25-cro1}, we indeed obtain that
\be\lab{2014-12-25-xe2}
 \lim_{n\rightarrow +\infty}\int_{\R^N}a_{i,\varepsilon_n}(x)|u_n|^{2^*(s_i)}dx=0\;\hbox{for all}\;i=1,2,\cdots,l.
\ee
Then by H\"older inequality and Hardy-Sobolev inequality, we see that
\be\lab{2014-12-17-we10}
\int_{\R^N}\Big[\sum_{i=1}^{l}\lambda_ia_{i,\varepsilon_n}(x)|u_n|^{2^*(s_i)-2}u_n h\Big]dx=o(1)\|h\|.
\ee
Recalling that $\Phi'_{\varepsilon_n}(u_n)=0$, we obtain that
\be\lab{2014-12-17-we11}
\langle \Psi'(u_n), h\rangle \equiv \int_{\R^N}\Big[\sum_{i=1}^{l}\lambda_ia_{i,\varepsilon_n}(x)|u_n|^{2^*(s_i)-2}u_n h\Big]dx=o(1)\|h\|.
\ee
Hence, $\Psi'(u_n)\rightarrow 0$.
\ep

\bc\lab{2014-12-17-cro1}
$\displaystyle \lim_{n\rightarrow +\infty}\int_{\R^N}\Big[a_{i,\varepsilon_n}(x)|u_n|^{2^*(s_i)}\Big]dx>0$ for all $i=1,2,\cdots,l$.
\ec
\bp
We prove it by the way of negation. We assume that
\be\lab{assumption}
\lim_{n\rightarrow +\infty}\int_{\R^N}\Big[a_{i,\varepsilon_n}(x)|u_n|^{2^*(s_i)}\Big]dx=0\;\hbox{for some}\;i\in \{1,2,\cdots,l\}.
\ee
Then by Lemma \ref{2014-12-17-l3}, $\{u_n\}$ is a $PS$ sequence of $\Psi$. By Remark \ref{2014-12-16-r1}, we always have $\displaystyle\liminf_{n\rightarrow +\infty}c_{\varepsilon_n}>0$.
Hence, $u_n\not\rightarrow 0$ in $D_{0}^{1,2}(\R^N)$, and then it is easy to see that
\be\lab{2014-12-17-we12}
\lim_{n\rightarrow+\infty}\Psi(u_n)\geq \frac{1}{N}S^{\frac{N}{2}}.
\ee
We note that under the assumption \eqref{assumption}, one can easily obtain that
\be\lab{2014-12-17-we13}
\lim_{n\rightarrow +\infty}\int_{\R^N}\Big[\sum_{i=1}^{l}\frac{1}{2^*(s_i)}\lambda_ia_{i,\varepsilon_n}(x)|u_n|^{2^*(s_i)}\Big]dx=0.
\ee
Thus, up to a subsequence, we can obtain that
\be\lab{2014-12-17-we14}
\lim_{n\rightarrow +\infty}\Phi_{\varepsilon_n}(u_n)=\lim_{n\rightarrow +\infty}\bigg[\Psi(u_n)-\int_{\R^N}\Big[\sum_{i=1}^{l}\frac{1}{2^*(s_i)}\lambda_ia_{i,\varepsilon_n}(x)|u_n|^{2^*(s_i)}\Big]dx\bigg]
\ee
and the above limit  is $\geq \frac{1}{N}S^{\frac{N}{2}},$ a contradiction to Lemma \ref{2014-12-16-l6}.
\ep

\bl\lab{2014-12-18-l1}
Let $\varepsilon_n\downarrow 0$.
Assume that $\{\phi_n\}\subset D_{0}^{1,2}(\R^N)$ is a bounded sequence such that
\be\lab{2014-12-18-e1}
\lim_{n\rightarrow +\infty}J_{\varepsilon_n}(\phi_n)=0
\ee
and $\phi_n\not\rightarrow 0$ in $L^{2^*}(\R^N)$. Suppose that  there exists some $i\in \{1,2,\cdots,l\}$ such that
\be\lab{2014-12-18-e2}
\liminf_{n\rightarrow +\infty}\int_{\R^N}\Big[a_{i,\varepsilon_n}(x)|\phi_n|^{2^*(s_i)}\Big]dx>0,
\ee
Then up to a subsequence, we must have
\be\lab{2014-12-18-e3}
\lim_{n\rightarrow +\infty}\Phi_{\varepsilon_n}(\phi_n)\geq \lim_{n\rightarrow +\infty} c_{\varepsilon_n}>0.
\ee
\el
\bp
Up to a subsequence if necessary, we denote
\be\lab{2014-12-17-we15}
\eta_i:=\liminf_{n\rightarrow +\infty} \int_{\R^N}\lambda_i a_{i,\varepsilon_n}(x)|\phi_n|^{2^*(s_i)}dx.
\ee
Obviously, $\phi_n\not\rightarrow 0$ in $D_{0}^{1,2}(\R^N)$. If not, by the Sobolev inequality we obtain that $\phi_n\rightarrow 0$ in $L^{2^*}(\R^n)$, a contradiction.
Since also that $\{\phi_n\}$ is bounded in $D_{0}^{1,2}(\R^N)$,  up to a subsequence,  there exists some $d_1,d_2>0$ such that
\be\lab{2014-12-17-we16}
0<d_1\leq \|\phi_n\|^2\leq d_2.
\ee
By the way, $\phi_n\not\rightarrow 0$ in $L^{2^*}(\R^N)$ yields that
there exist some $d_3>0$ such that
\be\lab{2014-12-17-we17}
d_3\leq |\phi_n|_{2^*}^{2^*}.
\ee
On the other hand, by the Sobolev inequality again, there exists some $d_4>0$ such that
\be\lab{2014-12-17-we20}
|\phi_n|_{2^*}^{2^*}\leq d_4.
\ee
Now, up to a subsequence, we may assume that
\be\lab{2014-12-18-e5}
\|\phi_n\|^2\rightarrow a^*>0, |\phi_n|_{2^*}^{2^*}\rightarrow b^*>0.
\ee
Then by the assumption \eqref{2014-12-18-e1}, we have that
\be\lab{2014-12-18-e6}
a^*-\sum_{i=1}^{l}\eta_i-b^*=0.
\ee
If there exists some $i\in \{1,2,\cdots,l\}$ such that \eqref{2014-12-18-e2} holds, then by Corollary \ref{2014-12-25-cro1}, we obtain that \eqref{2014-12-18-e2} holds for all $i\in \{1,2,\cdots,l\}$.
On the other hand, by Lemma \ref{2014-12-16-l1}, for $\phi_n$, there exists a  unique $t_n>0$ such that $t_n\phi_n\in \mathcal{N}_{\varepsilon_n}$.
Hence,
\be\lab{2014-12-18-e7}
\|\phi_n\|^2-\sum_{i=1}^{l}\lambda_ia_{i,\varepsilon_n}(x)|\phi_n|^{2^*(s_i)}t_{n}^{2^*(s_i)-2}-|\phi_n|_{2^*}^{2^*}t_{n}^{2^*-2}=0.
\ee
Then firstly we have
\begin{align}\lab{2014-12-18-e8}
\|\phi_n\|^2=&\sum_{i=1}^{l}\lambda_ia_{i,\varepsilon_n}(x)|\phi_n|^{2^*(s_i)}t_{n}^{2^*(s_i)-2}+|\phi_n|_{2^*}^{2^*}t_{n}^{2^*-2}\nonumber\\
\leq&C_i|\lambda_i|\|\phi_n\|^{2^*(s_i)}t_{n}^{2^*(s_i)-2}+C\|\phi_n\|^{2^*}t_{n}^{2^*-2}.
\end{align}
We claim that $t_n$ is bounded away from $0$. If not, we assume that $t_n\rightarrow 0$, then since $\|\phi_n\|\leq \sqrt{d_2}$, the right hand side of \eqref{2014-12-18-e8} goes to $0$. But by \eqref{2014-12-17-we16}, the left hand side of \eqref{2014-12-18-e8} is lager than $d_1>0$, we obtain a contradiction.
Secondly, by \eqref{2014-12-17-we17} and \eqref{2014-12-18-e7}, it is easy to see  that $\{t_n\}$ is bounded. Hence,   we may assume that $t_n\rightarrow t^*>0$.
Then we have

$$\left.
\begin{array}{lll}
&J_{\varepsilon_n}(t_n\phi_n)\equiv 0,\\
&\{\phi_n\}\;\hbox{is bounded in}\;D_{0}^{1,2}(\R^N),
\end{array}\right\}\Rightarrow \lim_{n\rightarrow +\infty}J_{\varepsilon_n}(t^*\phi_n)=0.
$$
Then it follows that
\be\lab{2014-12-18-e9}
a^*-\sum_{i=1}^{l}\eta_i (t^*)^{2^*(s_i)-2}-b^* (t^*)^{2^*-2}=0.
\ee
Apply the similar arguments of Lemma \ref{2014-12-16-l1},
we can  prove that the algebraic equation  $a^*-\sum_{i=1}^{l}\eta_i t^{2^*(s_i)-2}-b^* t^{2^*-2}=0$ has an unique positive solution. Hence, by \eqref{2014-12-18-e6} and \eqref{2014-12-18-e9}, we obtain that $t^*=1$.
Then by the boundedness of $\{\phi_n\}$ again, it is easy to see that
\be\lab{2014-12-18-e10}
\lim_{n\rightarrow +\infty}\Phi_{\varepsilon_n}(\phi_n)=\lim_{n\rightarrow +\infty}\Phi_{\varepsilon_n}(t_n\phi_n).
\ee
By the definition of $t_n$, we see that $t_n\phi_n\in \mathcal{N}_{\varepsilon_n}$. Hence, $\Phi_{\varepsilon_n}(t_n\phi_n)\geq c_{\varepsilon_n}$. It follows that
\be\lab{2014-12-18-e11}
\lim_{n\rightarrow +\infty}\Phi_{\varepsilon_n}(\phi_n)\geq \lim_{n\rightarrow +\infty}c_{\varepsilon_n}.
\ee
Insert  Remark \ref{2014-12-16-r1} here, we have that $\displaystyle \lim_{n\rightarrow +\infty}c_{\varepsilon_n}>0$.
\ep

\subsection{The proof of the existence result of Theorem \ref{2014-12-12-mainth2} for $k=l$}
Let $\varepsilon_n$ and $u_n$ be defined by Remark \ref{2015-2-16-r1}.
By Lemma \ref{2014-12-17-l2}, we only need to prove that $u_0\neq 0$. Now, we will proceed by contradiction. We assume that $u_0=0$.  By Corollary \ref{2014-12-17-cro1},
$$
\lim_{n\rightarrow +\infty}\int_{\R^N}\Big[a_{i,\varepsilon_n}(x)|u_n|^{2^*(s_i)}\Big]dx>0\;\hbox{for all}\;i=1,2,\cdots,l.
$$
Recalling that $\{u_n\}$ is bounded and all $\lambda_is$ are positive, up to a subsequence, we can denote that
\be\lab{2014-12-18-e12}
\lim_{n\rightarrow +\infty}\int_{\R^N}\Big[\sum_{i=1}^{l}\frac{\lambda_ia_{i,\varepsilon_n}(x)}{2^*(s_i)}|u_n|^{2^*(s_i)}\Big]dx=:\tau> 0.
\ee
Thus, by Corollary \ref{2014-12-16-cro2}, we obtain that
\begin{align}\lab{2014-12-18-e13}
&\lim_{n\rightarrow +\infty}\int_{\mathbb{B}_1}\Big[\sum_{i=1}^{l}\frac{\lambda_ia_{i,\varepsilon_n}(x)}{2^*(s_i)}|u_n|^{2^*(s_i)}\Big]dx\nonumber\\
&=\lim_{n\rightarrow +\infty}\int_{\mathbb{B}_1^c}\Big[\sum_{i=1}^{l}\frac{\lambda_ia_{i,\varepsilon_n}(x)}{2^*(s_i)}|u_n|^{2^*(s_i)}\Big]dx\nonumber\\
&=\frac{\tau}{2}> 0.
\end{align}
Let $\chi(x)\in C_c^\infty(\R^N)$ be a cut-off function such that $\chi(x)\equiv 1$ in $\mathbb{B}_{\frac{1}{2}}$, $\chi(x)\equiv 0$ in $\R^N\backslash \mathbb{B}_{1}$ and take $\tilde{\chi}(x)\in C^\infty(\R^N)$ such that $\tilde{\chi}(x)\equiv 0$ in $\mathbb{B}_1$ and $\tilde{\chi}\equiv 1$ in $\R^N\backslash \mathbb{B}_{2}$.
Let us denote
\be\lab{2014-12-18-e14}
\phi_{1,n}(x):=\chi(x)u_n(x), \phi_{2,n}(x):=\tilde{\chi}(x)u_n(x)
\ee
and
define
\be\lab{2014-12-18-e15}
\tilde{u}_n:=u_n-\phi_{1,n}-\phi_{2,n}.
\ee
Then we see that $sppt(\tilde{u}_n)\subset \Omega$, where $\Omega:=\{x\in \R^N:\;\frac{1}{2}<|x|<2\}$.
Then by the Rellich-Kondrachov compactness theorem, we see that
$\displaystyle \tilde{u}_n\rightarrow 0$ strongly in $L^{2^*(s_i)}(\Omega, \frac{dx}{|x|^{s_i}})$ for all $i=1,2,\cdots,l$.
Then it follows that $\tilde{u}_n$ is a $PS$ sequence of $\Psi$. By Br\'ezis-Lieb Lemma, we can  prove that
\be\lab{2014-12-18-e16}
\Phi_{\varepsilon_n}(u_n)=\Phi_{\varepsilon_n}(\phi_{1,n}+\phi_{2,n})+\Psi(\tilde{u}_n)+o(1).
\ee
Recalling that $\Phi'_{\varepsilon_n}(u_n)\equiv 0$, it is easy to prove that
\be\lab{2014-12-18-e17}
\lim_{n\rightarrow +\infty}\Phi'_{\varepsilon_n}(\phi_{1,n}+\phi_{2,n})=0.
\ee
Obviously,
\be\lab{2014-12-18-e18}
\lim_{n\rightarrow +\infty}\Phi_{\varepsilon_n}(\phi_{1,n}+\phi_{2,n})\geq 0.
\ee
Hence, if $\tilde{u}_n\not\rightarrow 0$ in $D_{0}^{1,2}(\R^N)$, we have that
\be\lab{2014-12-18-e19}
\lim_{n\rightarrow +\infty}\Psi(\tilde{u}_n)\geq \frac{1}{N}S^{\frac{N}{2}}.
\ee
By \eqref{2014-12-18-e17}, \eqref{2014-12-18-e18} and \eqref{2014-12-18-e16}, we obtain that $\displaystyle\lim_{n\rightarrow +\infty}c_{\varepsilon_n}\geq \frac{1}{N}S^{\frac{N}{2}}$, a contradiction to Lemma \ref{2014-12-16-l6}.
Hence, we prove that $\tilde{u}_n\rightarrow 0$ in $D_{0}^{1,2}(\R^N)$ and it follows that
\be\lab{2014-12-18-e20}
\Phi_{\varepsilon_n}(u_n)=\Phi_{\varepsilon_n}(\phi_{1,n})+\Phi_{\varepsilon_n}(\phi_{2,n})+o(1).
\ee
Recal that $\Phi'_{\varepsilon_n}(u_n)\equiv 0$ and hence
$
\langle \Phi'_{\varepsilon_n}(u_n), \phi_{1,n}\rangle \equiv 0.$
Then by $\tilde{u}_n\rightarrow 0$ strongly in $D_{0}^{1,2}(\R^N)$, it is easy to see that
\be\lab{2014-12-18-e22}
\lim_{n\rightarrow +\infty}J_{\varepsilon_n}(\phi_{1,n})=0.
\ee
By \eqref{2014-12-18-e13} and the Rellich-Kondrachov compactness result, we can prove that
\begin{align}\lab{2014-12-18-e23}
&\liminf_{n\rightarrow +\infty}\int_{\R^N}\Big[\sum_{i=1}^{l}\frac{\lambda_ia_{i,\varepsilon_n}(x)}{2^*(s_i)}|\phi_{1,n}|^{2^*(s_i)}\Big]dx\nonumber\\
&=\lim_{n\rightarrow +\infty}\int_{\mathbb{B}_1}\Big[\sum_{i=1}^{l}\frac{\lambda_ia_{i,\varepsilon_n}(x)}{2^*(s_i)}|u_n|^{2^*(s_i)}\Big]dx\nonumber\\
&=\frac{\tau}{2}> 0.
\end{align}
And it follows easily that
\be\lab{2014-12-26-e5}
\liminf_{n\rightarrow +\infty}\int_{\R^N}\Big[a_{i,\varepsilon_n}(x)|\phi_{1,n}|^{2^*(s_i)}\Big]dx>0\;\hbox{for all}\;i=1,2,\cdots,l.
\ee
Hence, by \eqref{2014-12-18-e22}, \eqref{2014-12-18-e23}, \eqref{2014-12-26-e5}  and Lemma \ref{2014-12-18-l1}, we obtain that
\be\lab{2014-12-18-e24}
\liminf_{n\rightarrow +\infty}\Phi_{\varepsilon_n}(\phi_{1,n})\geq \lim_{n\rightarrow +\infty}c_{\varepsilon_n}.
\ee
Similarly, we can also obtain that
\be\lab{2014-12-18-e25}
\liminf_{n\rightarrow +\infty}\Phi_{\varepsilon_n}(\phi_{2,n})\geq \lim_{n\rightarrow +\infty}c_{\varepsilon_n}.
\ee
Hence, by \eqref{2014-12-18-e20}, \eqref{2014-12-18-e24} and \eqref{2014-12-18-e25}, we have that
$\displaystyle
\lim_{n\rightarrow +\infty}c_{\varepsilon_n} \geq 2\lim_{n\rightarrow +\infty}c_{\varepsilon_n},
$ it is
a contradiction to the fact of that $\displaystyle\lim_{n\rightarrow +\infty}c_{\varepsilon_n}>0$  because of  Lemma \ref{2014-12-18-l1}. Thereby $u_0\neq 0$ is proved.\hfill$\Box$

\subsection{The proof of the existence result of Theorem \ref{2014-12-12-mainth2} for $k\neq l$}

 When $k\neq l$,  the proof becomes very thorny and we have to  apply another way-perturbation methods. In this case, we assume that $l\geq 2$.
 For the convenience, in this subsection  we  denote
\be\lab{2015-3-21-we1}
I_0(u)=I(u):=\frac{1}{2}\int_{\R^N}\left|\nabla u\right|^2 dx-\int_{\R^N}\sum_{i=1}^{k}\frac{\lambda_i}{2^*(s_i)}\frac{|u|^{2^*(s_i)}}{|x|^{s_i}},\;u\in D_{0}^{1,2}(\R^N),
\ee
and
\be\lab{2015-3-21-we2}
I_\lambda (u)=I_0(u)-\lambda \int_{\R^N}\left(\sum_{i=k+1}^{l}\frac{1}{2^*(s_i)} |\lambda_i| \int_{\R^N}\frac{|u|^{2^*(s_i)}}{|x|^{s_i}}\right),\;u\in D_{0}^{1,2}(\R^N),
\ee
which is the corresponding functional of the following variant problem:
\be\lab{2015-3-28-xe1}
\begin{cases}
\Delta u+u^{2^*-1}+\sum_{i=1}^{k}\lambda_i \frac{u^{2^*(s_i)-1}}{|x|^{s_i}}+\lambda\left(\sum_{i=k+1}^{l}|\lambda_i| \frac{u^{2^*(s_i)-1}}{|x|^{s_i}}\right)=0,\;x\in\;\R^N,\\
u\in D_{0}^{1,2}(\R^N).
\end{cases}
\ee
We note that when $\lambda=-1$, it becomes the problem \eqref{PMS}. Hence, our aim is to prove that problem \eqref{2015-3-28-xe1} possesses a least energy solution when $\lambda=-1$. We set
\be
D_k:=\left\{\mu\in \R\;:\;\hbox{problem \eqref{2015-3-28-xe1} possesses a least energy solution when}\;\lambda=\mu\right\}.
\ee
Then we only need to prove that $-1\in D_k$. Firstly, by the results established  in the  previous subsection, for any  $\lambda\geq 0$, problem \eqref{2015-3-28-xe1} possesses a least energy solution, and the corresponding energy is less than  $\displaystyle \frac{1}{N}S^{\frac{N}{2}}$.
Hence,
\be\lab{2014-3-21-rong1}
[0,+\infty)\subset D_k.
\ee
Set
\be\lab{2015-3-21-wmami1}
\mathcal{A}_{\mu}:=\left\{\hbox{$\displaystyle u\in D_{0}^{1,2}(\R^N)$ is a positive solution of problem \eqref{2015-3-28-xe1} when  $\lambda=\mu$}\right\}.
\ee
When $\displaystyle \mathcal{A}_{\lambda}\neq \emptyset$, we define
\be\lab{2015-3-21-wmami2}
c_{\lambda}^{*}:= \inf_{u\in \mathcal{A}_{\lambda}} I_{\lambda}(u),
\ee
\be\lab{2015-3-21-rong2}
c_\lambda:=\inf_{u\in D_{0}^{1,2}(\R^N)\backslash \{0\}}\max_{t>0}I_\lambda(tu)
\ee
and
\be\lab{2015-3-22-e2}
\tilde{c}_\lambda:=\min_{\gamma\in \Gamma_\lambda}\max_{t\in [0,1]}I_\lambda(\gamma(t)),
\ee
 where $\displaystyle\Gamma_\lambda:=\left\{\gamma\in C\left([0,1], D_{0}^{1,2}(\R^N)\right):   \; \gamma(0)=0, I_\lambda(\gamma(1))<0\right\}$.
Then it is easy to see that
\be\lab{2015-3-23-lala5}
 c_{\lambda}^{*}\geq c_{\lambda}\geq \tilde{c}_{\lambda}>0.
 \ee
 By the standard concentration compactness arguments, one can prove that if there exists a bounded $(PS)_c$ sequence of $I_\lambda$ with $c<c_{\lambda}^{*}$, then $c$ is a critical value of $I_\lambda$. It follows that
\be\lab{2015-3-23-lala8}
 c_{\lambda}^{*}= c_{\lambda}= \tilde{c}_{\lambda}>0.
 \ee
Next,  we prepare  the following properties about the least energy solution.
\bl\lab{2015-3-21-l1}
Assume that for some  $\lambda\in \R$, equation  \eqref{2015-3-28-xe1} possesses a least energy solution $u_\lambda(x)$, then
\be\lab{2015-3-21-we3}
0<I_\lambda(u_\lambda)=c_\lambda <\frac{1}{N}S^{\frac{N}{2}}.
\ee
On the other hand,  the functional $I_\lambda$  possesses  the following properties:
\begin{itemize}
\item[(M1)] there exists some   $c, r>0$ such that  $I_\lambda(u)\geq c$ for   $\|u\|=r$.   Moreover, there exists  $v_\lambda\in D_{0}^{1,2}(\R^N)$  such that  $\|v_\lambda\|>r$ and $I_\lambda(v_\lambda)<0$;
\item[(M2)] there exists a critical point $u_\lambda\in D_{0}^{1,2}(\R^N)$ of $I_\lambda$ such that
    \be I_\lambda(u_\lambda)=c_\lambda=\tilde{c}_\lambda:=\min_{\gamma\in \Gamma}\max_{t\in [0,1]}I_\lambda(\gamma(t)),\ee
     where $\displaystyle\Gamma_\lambda:=\left\{\gamma\in C\left([0,1], D_{0}^{1,2}(\R^N)\right):   \; \gamma(0)=0, \gamma(1)=v_\lambda\right\}$;
\item[(M3)]  $ c_\lambda=c_{\lambda}^{*}=\inf\left\{I_\lambda(u):  \; I'_\lambda(u)=0, u\in D_{0}^{1,2}(\R^N)\backslash \{0\}\right\}$；
\item[(M4)] the set $\mathcal{S}_\lambda:=\left\{u\in D_{0}^{1,2}(\R^N):    \; I'_\lambda(u)=0, I_\lambda(u)=c_\lambda\right\}$ is compact in $D_{0}^{1,2}(\R^N)$;
\item[(M5)]there exists a path $\gamma_\lambda(t)\in \Gamma_\lambda$ passing through $u_\lambda$ at $t=t_\lambda$ and satisfying
    \be I_\lambda(u_\lambda)>I_\lambda\left(\gamma_\lambda(t)\right)\;\hbox{for all $t\neq t_\lambda$}.\ee
\end{itemize}
\el
\bp Obviously, $c_\lambda>0$. Combining with the result of Lemma \ref{2014-12-16-l6}, we obtain  \eqref{2015-3-21-we3} and \eqref{2015-3-21-rong2}. Based  on the result of \eqref{2015-3-23-lala8}, (M1)-(M3) are trivial. And by Lemma \ref{2014-12-16-l1} we can obtain  (M5). Hence, next we only need to check the property of (M4). Let $\{u_n\}\subset \mathcal{S}_\lambda$, noting that $I'_\lambda(u_n)=0$, by Lemma \ref{2014-12-16-l2} we see that $\{u_n\}$ is a bounded $(PS)_{c_\lambda}$ sequence of $I_\lambda$. And it is easy to prove that $|u_n|_{2^*}$ are bounded away from $0$.  By Proposition \ref{2014-12-20-prop1}, we see that $\{u_n(0)\}$ is bounded. Hence,  $u_n(x)$ is a bounded sequence of $L^\infty(\R^N)\cap D_{0}^{1,2}(\R^N)$. Noting that the Kelvin transform of $u_n$, which is denoted by $\hat{u}_n(x):=|x|^{-(N-2)}u_n\left(\frac{x}{|x|^2}\right)$, is also a least energy solution, i.e.,  $\hat{u}_n\in \mathcal{S}_\lambda$. On the other hand, we also note that for any $s\in [0,2]$  and any solution $u$ with its Kelvin transform $\hat{u}$, we have
\be\lab{2015-3-22-ze1}
\int_{\mathbb{B}_1} \frac{|u|^{2^*(s)}}{|x|^s}dx=\int_{\mathbb{B}_1^c} \frac{|\hat{u}|^{2^*(s)}}{|x|^s}dx.
\ee
Hence, for the new sequence $\displaystyle \{u_1, \hat{u}_1, u_2, \hat{u}_2, \cdots\}\subset \mathcal{S}_\lambda$, there exists a subsequence, denoted by  $w_j$, such that
\be\lab{2015-3-22-ze2}
\liminf_{j\rightarrow \infty}\int_{\mathbb{B}_1} |w_j|^{2^*}dx>0.
\ee
Then by Lemma \ref{2014-12-18-wl1}, we obtain that
\be\lab{2015-3-22-ze3}
\liminf_{j\rightarrow +\infty} w_j(0)>0,
\ee
Hence, up to a subsequence,  $w_j \rightharpoonup w\neq 0$ in $D_{0}^{1,2}(\R^N)$. It is easy to see that $w$ is also a critical point of  $I_\lambda$. Hence, we have
\be\lab{2015-3-22-ze4}
I_\lambda(w)\geq c_\lambda.
\ee
On the other hand, by the weak semi-continuous of a norm, when $\lambda\leq 0$,
\begin{align}\lab{2015-3-22-ze5}
I_\lambda(w)=& \left[\frac{1}{2}-\frac{1}{2^*(s_k)}\right]\|w\|^2+\sum_{i=1}^{k-1}\left[\frac{1}{2^*(s_k)}-\frac{1}{2^*(s_i)}\right]
\lambda_i|w|_{2^*(s_i)}^{2^*(s_i)}\nonumber\\
&+\sum_{i=k+1}^{l}\lambda \left[\frac{1}{2^*(s_k)}-\frac{1}{2^*(s_i)}\right]
|\lambda_i||w|_{2^*(s_i)}^{2^*(s_i)}+\left[\frac{1}{2^*(s_k)}-\frac{1}{2^*}\right]|w|_{2^*}^{2^*}\\
\leq&\liminf_{j\rightarrow \infty}\left\{\left[\frac{1}{2}-\frac{1}{2^*(s_k)}\right]\|w_j\|^2+\sum_{i=1}^{k-1}\left[\frac{1}{2^*(s_k)}-\frac{1}{2^*(s_i)}\right]
\lambda_i|w_j|_{2^*(s_i)}^{2^*(s_i)}\right.\nonumber\\
&\left.+\sum_{i=k+1}^{l}\lambda \left[\frac{1}{2^*(s_k)}-\frac{1}{2^*(s_i)}\right]
|\lambda_i||w_j|_{2^*(s_i)}^{2^*(s_i)}+\left[\frac{1}{2^*(s_k)}-\frac{1}{2^*}\right]|w_j|_{2^*}^{2^*}\right\}\nonumber\\
=&\liminf_{i\rightarrow \infty} I_{\lambda}(w_j)=c_\lambda.\nonumber
\end{align}
The case of $\lambda>0$ is much easier to check. Then, it follows that $I_\lambda(w)=c_\lambda$ and thus $w_j\rightarrow w\in \mathcal{S}_\lambda$, its Kelvin transform $\hat{w}$ also satisfies  $\hat{w}\in \mathcal{S}_\lambda$. Hence, up to a subsequence, we have $u_n\rightarrow w\in \mathcal{S}_\lambda$ or $u_n\rightarrow \hat{w}\in \mathcal{S}_\lambda$. Thereby, (M4) is verified.
\ep


\bl\lab{2015-3-21-l2}
For any $\lambda\in D_k$, there exists some $\delta_\lambda>0$ small enough such that $\displaystyle (\lambda-\delta_\lambda, \lambda+\delta_\lambda)\subset D_k$. In other words, $D_k$ is an open set of $\R$.
\el
\bp
Basing on the Lemma \ref{2015-3-21-l1}, applying the perturbation arguments, it is standard to prove the existence of $\delta_\lambda$ and the existence of positive solution for $\displaystyle\mu\in (\lambda-\delta_\lambda, \lambda+\delta_\lambda)$. This processes is very long and tedious, however it is standard. Hence, we omit the details and a very like discussion we refer to \cite[section 5]{CeramiZhongZou.2015}. Next, we shall prove the existence of least energy solution. For any fixed $\mu\in (\lambda-\delta_\lambda, \lambda+\delta_\lambda) $, then we firstly have that $\displaystyle \mathcal{A}_{\mu}\neq \emptyset$.
Let  $\{u_n\}\subset \mathcal{A}_{\mu}$ be a minimizing sequence, then it is easy to see that $\{u_n\}$ is a bounded  $(PS)_{c_{\mu}^{*}}$ sequence of $I_{\mu}$. Let  $\hat{u}_n$ be the Kelvin transform of $u_n$, then we also have that $\{\hat{u}_n\}\subset \mathcal{A}_{\mu}$. Hence, apply the similar argument of the   (M4) in Lemma \ref{2015-3-21-l1}, we can prove that $\displaystyle \{u_1, \hat{u}_1, u_2, \hat{u}_2, \cdots\}$ is also a minimizing sequence and it possesses a strong convergent subsequence. Thus, we prove that $c_{\mu}^{*}$ is achieved. We also note that \eqref{2015-3-23-lala8} holds.
\ep

\bl\lab{2015-3-21-l3}
$D_k$ is closed in $\R$.
\el
\bp
For any sequence $\{\lambda_n\}\subset D_k, \lambda_n\rightarrow \lambda$, we shall prove that $\lambda\in D_k$.

Firstly we have
\be
\lim_{n\rightarrow \infty}c_{\lambda_n}=\lim_{n\rightarrow \infty}\tilde{c}_{\lambda_n}= \tilde{c}_{\lambda}\leq c_\lambda.
\ee
Secondly by Lemma \ref{2014-12-16-l6},
\be
0<\lim_{n\rightarrow \infty}c_{\lambda_n}<\frac{1}{N}S^{\frac{N}{2}}.
\ee
Then the boundedness of $c_{\lambda_n}$ yields the boundedness of $\{u_n\}$ in $D_{0}^{1,2}(\R^N)$. By $\lambda_n\rightarrow \lambda$ and $I'_{\lambda_n}(u_n)=0$, we see that $I'_\lambda(u_n)\rightarrow 0$. Hence, we obtain that  $\{u_n\}$ is a bounded $(PS)_{\tilde{c}_\lambda}$ sequence of $I_\lambda$. We still adopt the notation  $\hat{u}_n$ as the Kelvin transform of $u_n$, then we firstly have that $I'_{\lambda_n}(\hat{u}_n)=0$, furthermore, we have that
$\{\displaystyle u_1, \hat{u}_1, u_2,\hat{u}_2,\cdots\}$ is also a bounded $(PS)_{\tilde{c}_\lambda}$ sequence of $I_\lambda$.
Noting that
$$\int_{\mathbb{B}_1} \frac{|u_i|^{2^*(s)}}{|x|^s}dx=\int_{\mathbb{B}_1^c} \frac{|\hat{u_i}|^{2^*(s)}}{|x|^s}dx\;\hbox{for all $i=1,2,\cdots$ and $s\in [0,2]$}.$$
Then applying the similar argument of the  (M4) in Lemma \ref{2015-3-21-l1}, we obtain that $\{\displaystyle u_1, \hat{u}_1, u_2,\hat{u}_2,\cdots\}$ possesses a strong convergent subsequence. Hence, up to a subsequence, we may assume that $u_n\rightarrow u$ or $u_n\rightarrow \hat{u}$. Hence, $\tilde{c}_\lambda$ is achievable, and it follows that $\mathcal{A}_{\lambda}\neq \emptyset$. Then we have the relationship of  \eqref{2015-3-23-lala8} and thus $u$ is a least energy solution.
Hence, $\lambda\in D_k$, and $D_k$ is closed in $\R$.
\ep

\vskip 0.2in
\noindent{\bf The  final proof of the existence of least energy solution of Theorem \ref{2014-12-12-mainth2} for $k\neq l$:}
By Lemma \ref{2015-3-21-l2} and Lemma \ref{2015-3-21-l3}, we see that $D_k$ is a both open and closed set of $\R$. By \eqref{2014-3-21-rong1},  $[0,+\infty)\subset D_k$, hence $D_k\neq \emptyset$. Finally,  we obtain that $D_k=\R$. Then $-1\in D_k$, and thus the existence of Theorem \ref{2014-12-12-mainth2} is completed.\hfill$\Box$

\vskip0.26in


\end{document}